\documentclass[12pt]{article}
\usepackage{latexsym, amssymb, amsmath, amscd, amsfonts, epsfig, graphicx, colordvi,verbatim,ifpdf}
\usepackage{amsfonts, amsmath, amssymb,extarrows}
\usepackage{amssymb,amsfonts,amsmath,latexsym,epsfig,cite, psfrag,eepic,color}
\usepackage{amscd,graphics}
\usepackage{latexsym, amssymb,  amsmath,amscd, amsfonts, epsfig, graphicx, colordvi,amsthm}

\usepackage{graphicx}
\usepackage{color}
\usepackage{ifpdf}
\usepackage{fancybox}

\usepackage{float}

\newtheorem{thm}{Theorem}[section]

\newtheorem{conj}[thm]{Conjecture}

\newtheorem{lem}[thm]{Lemma}
\newtheorem{core}[thm]{Corollary}

\numberwithin{equation}{section}

\allowdisplaybreaks

\newlength{\boxedparwidth}
  {\begin{center} \begin{tabular}{|@{\hspace{.315in}}c@{\hspace{.15in}}|}
                  \hline \\ \begin{minipage}[t]{\boxedparwidth}
                  \setlength{\parindent}{.25in}}%
  {\end{minipage} \\ \\ \hline \end{tabular} \end{center}}

\begin{document}

\begin{center}

 {\large \bf Unimodality of the Andrews-Garvan-Dyson cranks of  partitions}
\end{center}

\begin{center}
{Kathy Q. Ji}$^{1}$ and
  {Wenston J.T. Zang}$^{2}$ \vskip 2mm

  $^{1}$Center for Applied Mathematics,\\[2pt]
Tianjin University, Tianjin 300072, P.R. China\\[3pt]

   \vskip 2mm

   $^{2}$Institute for Advanced Study in Mathematics,\\[2pt] Harbin Institute of Technology, Heilongjiang, 150001, P.R. China\\[8pt]
  $^1$kathyji@tju.edu.cn, \quad $^2$zang@hit.edu.cn
\end{center}

\vskip 6mm \noindent {\bf Abstract.} The main objective of this paper is to investigate the
distribution of the Andrews-Garvan-Dyson cranks of  partitions. Let
$M(m,n)$ denote the number of partitions of $n$ with the Andrews-Garvan-Dyson
crank $m$, we show that the sequence
$\{M(m,n)\}_{|m|\leq n-1}$ is unimodal for $n\geq 44$. It turns out that
the unimodality of
$\{M(m,n)\}_{|m|\leq n-1}$ is related to the monotonicity properties of
two partition functions $p_{r}(n)$ and $pp_{r}(n)$. Let $p_{r}(n)$ denote
the number of partitions of $n$ with parts taken from
$\{2,3,\ldots ,r\}$ and let $pp_{r}(n)$ denote the number of pairs
$(\alpha ,\beta )$ of partitions, where $\alpha $ is a partition counted
by $p_{r}(i)$ and $\beta $ is a partition counted by $p_{r+1}(n-i)$ for
$0\leq i\leq n$. We show that $p_{r}(n)\geq p_{r}(n-1)$ for
$r\geq 5$ and $n\geq 14$ and $pp_{r}(n)\geq pp_{r}(n-1)$ for
$r\geq 3$ and $n\geq 8$. With the aid of the monotonicity properties on
$p_{r}(n)$ and $pp_{r}(n)$, we show that $M(m,n)\geq M(m,n-1)$ for
$n\geq 14$ and $ 0\leq m \leq n-2$ and $M(m-1,n)\geq M(m,n)$ for
$n\geq 44$ and $1\leq m\leq n-1$. By means of the symmetry
$M(m,n)=M(-m,n)$, we find that $M(m-1,n)\geq M(m,n)$ for $n\geq 44$ and
$1\leq m\leq n-1$ implies that the sequence
$\{M(m,n)\}_{|m|\leq n-1}$ is unimodal for $n\geq 44$. We also give a proof
of an upper bound for $\mathop{\mathrm{ospt}}\nolimits (n)$ conjectured
by Chan and Mao in light of $M(m-1,n)\geq M(m,n)$ for $n\geq 44$ and
$0\leq m\leq n-1$.
$\\[5pt]$
{{\bf Keywords}: Partition, rank, crank, unimodal, ospt-function}\\[3pt]
{\bf 2010 Mathematics Subject Classification.} {11P81, 05A17, 05A20}

 \section{Introduction}
\label{sec1}

Dyson's rank \cite{Dyson-1944} and the Andrews-Garvan-Dyson crank
\cite{Andrews-Garvan-1988} are two fundamental statistics in the theory
of partitions. Recall that the rank of a partition was introduced by Dyson
\cite{Dyson-1944} as the largest part of the partition minus the number
of parts. The crank of a partition was defined by Andrews and Garvan
\cite{Andrews-Garvan-1988} as the largest part if the partition contains
no ones, and otherwise as the number of parts larger than the number of
ones minus the number of ones.

Let $p(n)$ denote the number of partitions of $n$. It was conjectured by
Dyson \cite{Dyson-1944} and confirmed by Atkin and Swinnerton-Dyer
\cite{Atkin-Swinnerton-Dyer-1954} that the rank of a partition could explain
two of Ramanujan's famous partition congruences
$p(5n+4)\equiv 0 \pmod{5}$ and $p(7n+5)\equiv 0 \pmod{7}$, but not the
third one, $p(11n+6)\equiv 0\pmod{11}$. This led Dyson to hypothesize the
existence of another statistic, namely the crank. Until forty-four years
later, Andrews and Garvan \cite{Andrews-Garvan-1988}, building on the work
of Garvan \cite{Garvan-1988} finally unveiled crank and showed that the
crank can be used to interpret all three congruences on $p(n)$ mod
$5,7$ and $11$. For more details, please refer to Dyson
\cite{Dyson-1944}, Atkin and Swinnerton-Dyer
\cite{Atkin-Swinnerton-Dyer-1954} and Andrews and Garvan
\cite{Andrews-Garvan-1988, Garvan-1988}. It is worth mentioning that Mahlburg
\cite{Mahlburg-2005} showed that the crank can also provide combinatorial
interpretations of infinite families of congruences on $p(n)$ established
by Ahlgren and Ono \cite{Ahlgren-Ono-2001} and Ono \cite{Ono-2000}. Since
then, the rank and the crank have been extensively studied, see, for example,
Andrews and Garvan \cite{Andrews-1989-242}, Andrews and Ono
\cite{Andrews-Ono-2005}, Bringmann and Dousse
\cite{Bringmann-Dousse-2016}, Bringmann and Ono
\cite{Bringmann-Ono-2006, Bringmann-Ono-2010}, Garvan
\cite{Garvan-1990}, Lewis \cite{Lewis-1991b}, and so on.

Let $m$ be an integer. For $n\geq 1$, let $N(m,n)$ denote the number of
partitions of $n$ with rank $m$, and for $n>1$, let $M(m,n)$ denote the
number of partitions of $n$ with crank $m$. For $n=1$, set
$M(0,1)=-1,\,M(1,1)=M(-1,1)=1$, and $ M(m,1)=0 $ when
$m\neq -1,0,1$. For $n=0$, set $M(0,0)=1$, and $ M(m,0)=0 $ when
$m\neq 0$. For $n<0$, set $M(m,n)=0$.

In 2014, Chan and Mao \cite{Chan-Mao-2014} showed the following two inequalities
on $N(m,n)$:
%
\begin{thm}{\rm(Chan and Mao).}%
\label{thm-cm-1}
For $n\geq 12$ and $ 0\leq m\leq n-3$ or $m=n-1$,
%
\begin{equation}
N(m,n)\geq N(m,n-1).
\end{equation}
\end{thm}

\begin{thm}{\rm(Chan and Mao).}%
\label{thm-cm-2}
For $m,n\geq 0$,
%
\begin{equation}
N(m,n)\geq N(m+2,n).
\end{equation}
\end{thm}

In \cite{Andrews-Chan-Kim-2013}, Andrews, Chan and Kim introduced the function
$\mathop{\mathrm{ospt}}\nolimits (n)$ defined as the difference between
the first positive crank moment and the first positive rank moment, namely,
%
\begin{equation}
\label{equ-ospt}
\mathop{\mathrm{ospt}}\nolimits (n)=\sum _{m=0}^{\infty }mM(m,n)-
\sum _{m=0}^{\infty }mN(m,n).
\end{equation}
By means of generating function, Andrews, Chan and Kim
\cite{Andrews-Chan-Kim-2013} proved the positivity of
$\mathop{\mathrm{ospt}}\nolimits (n)$ and gave a combinatorial interpretation
of $\mathop{\mathrm{ospt}}\nolimits (n)$ which counts the number of even
and odd strings in the partitions of $n$. Chen, Ji and Zang
\cite{Chen-Ji-Zang-2017} gave another combinatorial interpretation of
$\mathop{\mathrm{ospt}}\nolimits (n)$ in terms of certain bijection.

Using  {Theorem~\ref{thm-cm-1}} and  {Theorem~\ref{thm-cm-2}}, Chan and Mao
\cite{Chan-Mao-2014} established the following upper-bound and lower-bound
for $\mathop{\mathrm{ospt}}\nolimits (n)$ in terms of $N(m,n)$,
$M(m,n)$ and $p(n)$.

\begin{thm}{\rm(Chan and Mao).}
The following inequalities are true.
%
\begin{align}
\mathop{\mathrm{ospt}}\nolimits (n)&>\frac{p(n)}{4}+
\frac{N(0,n)}{2}-\frac{M(0,n)}{4},&\text{for }n\geq 8,
\\[3pt]
\mathop{\mathrm{ospt}}\nolimits (n)&<\frac{p(n)}{4}+
\frac{N(0,n)}{2}-\frac{M(0,n)}{4}+\frac{N(1,n)}{2},&\text{for }n\geq 7,
\label{chan-mao-upbound}
\\[3pt]
\mathop{\mathrm{ospt}}\nolimits (n)&<\frac{p(n)}{2},&\text{for }n
\geq 3.
\end{align}
\end{thm}

At the end of the paper, Chan and Mao \cite{Chan-Mao-2014} raised a series
of open problems, one of which is to establish similar inequalities for
the crank of a partition. They also posed the following conjecture.

\begin{conj}{\rm(Chan and Mao).}%
\label{conj-opst3}
For $n\geq 10$,
%
\begin{equation}
\label{opstpn3}
\mathop{\mathrm{ospt}}\nolimits (n)<\frac{p(n)}{3}.
\end{equation}
\end{conj}

In \cite{Kim-Kim-Seo}, Kim, Kim and Seo proved that
$M(m,n)>M(m+1,n)$ for $m\geq 0$ and sufficiently large $n$. More precisely,
they obtained the following result.

\begin{thm}{\rm(Kim, Kim and Seo).}%
\label{thm-kks}
For $m\geq 0$,
\begin{equation*}
M(m,n)>M(m+1,n)
\end{equation*}
for all positive integers $n>100$ satisfying
\begin{equation*}
\sqrt{n}I_{-9/2}\left (\pi \sqrt{\frac{2n}{3}}\right )>217
\frac{(2m+3)^{14}}{2m+1}e^{\frac{\pi \sqrt{3}}{32} (2m+3)^{2}}e^{\pi
\sqrt{\frac{2n}{3}}},
\end{equation*}
where $I_{s}(z)$ is the modified Bessel function of the second kind.
\end{thm}
In this paper, we establish the following two inequalities on
$M(m,n)$.

\begin{thm}%
\label{main-un-2-n}
For $n\geq 14$ and $0\leq m\leq n-2$,
%
\begin{equation}
\label{equ-main-thm-n-un}
M(m,n)\geq M(m,n-1).
\end{equation}
\end{thm}

\begin{thm}%
\label{main-thm-n}
For $n\geq 44$ and $1\leq m \leq n-1$,
%
\begin{equation}
\label{main-ine}
M(m-1,n)\geq M(m,n).
\end{equation}
\end{thm}

Recall that a sequence $\{a_{i}\}_{1\le i\le n}$ is unimodal if for some
$1\le j\le n$,
\begin{equation*}
a_{1}\le \cdots \le a_{j-1}\le a_{j}\ge a_{j+1}\ge \cdots \ge a_{n}.
\end{equation*}
For more information, see \cite[P.124, Ex.50]{Stanley-1997}.

From  {Theorem~\ref{main-thm-n}} and the symmetry $M(m,n)=M(-m,n)$ (see
\cite{Garvan-1988,Dyson-1989}), we find the following unimodality of the
crank.

\begin{core}
\label{unimod-cor}
For $n\geq 44$,
\begin{equation*}
M(1-n,n)\leq \cdots \le M(-1,n) \leq M(0,n)\geq M(1,n)\geq \cdots
\geq M(n-1,n).
\end{equation*}
That means the sequence $\{M(m,n)\}_{|m|\leq n-1}$ is unimodal for
$n\geq 44$.
\end{core}

\begin{figure}[h]\label{unimodal}
\centering
\includegraphics[width=4in]{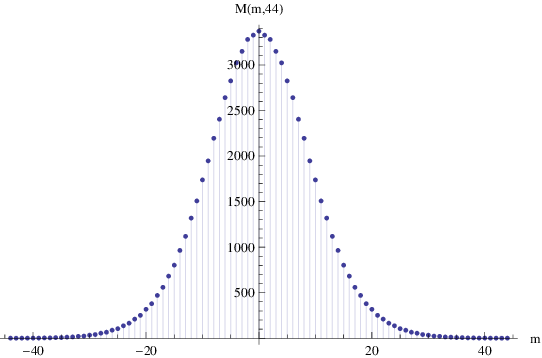}
\caption{The sequence $\{M(m,44)\}_{|m|\leq 43}$ is unimodal.}
\end{figure}

 {Fig. 1.1} gives an illustration of the unimodality of
$\{M(m,44)\}_{|m|\leq 43}$.

It should be noted that when
\begin{equation*}
(m,n)\in \{(2,5),(3,10),(4,9),(6,13)\},
\end{equation*}
the inequality  {\eqref{equ-main-thm-n-un}} does not hold.

Also, when
\begin{align*}
(m,n)\in &\{(1,2i-1)\colon 4\le i\le 22\}\cup \{(2,4)\}\cup \{(2,2i)
\colon 5\le i\le 13\}\cup
\\[3pt]
&\{(3,5),(3,9),(3,15),(3,17),(3,21),(4,8),(4,10),(4,16),(5,9),(5,13),(6,12)
\},
\end{align*}
the inequality  {\eqref{main-ine}} does not hold.

It is worth mentioning that Andrews, Dyson and Rhoades
\cite{Andrews-Dyson-Rhoades-2013} conjectured the unimodality of the spt-crank
defined on the spt-function. Let $N_{S}(m,n)$ denote the number of
$S$-partitions of $n$ with spt-crank $m$, Andrews, Dyson and Rhoades conjectured
that $\{N_{S}(m,n)\}_{m}$ is unimodal. Their conjecture was proved by Chen,
Ji and Zang \cite{Chen-Ji-Zang-2015}. For the definitions of the spt-crank,
the spt-function and the $S$-partition, please refer to
\cite{Andrews-2008} and \cite{Andrews-Garvan-Liang-2011}.

In this paper, we also give a proof of  {Conjecture~\ref{conj-opst3}} in light
of  {Theorem~\ref{main-thm-n}} and the following theorem.

\begin{thm}%
\label{lem-indu}
For $n\geq 39$,
%
\begin{equation}
\label{esm0n}
p(n)\geq 21M(0,n).
\end{equation}
\end{thm}

The proofs of the two inequalities on $M(m,n)$ relates to the monotonicity
property of the partition function $p_{r}(n)$, which counts the number
of partitions of $n$ with parts taken from $\{2,3,\ldots ,r\}$. From the
definition of $p_{r}(n)$, it is easy to see that the generating function
of $p_{r}(n)$ is
%
\begin{equation}
\label{gf-pkn}
\sum _{n=0}^{\infty }p_{r}(n)q^{n}=\frac{1}{(q^{2};q)_{r-1}}.
\end{equation}
Here and throughout the rest of this paper, we adopt the common $q$-series
notation \cite{and76}:
\begin{align*}
(a;q)_{\infty }=\prod _{n=0}^{\infty }(1-aq^{n}) \quad \text{and} \quad (a;q)_{n}&=
\frac{(a;q)_{\infty }}{(aq^{n};q)_{\infty }}.
\end{align*}

We show that $p_{r}(n)$ has the following monotonicity property.

\begin{thm}%
\label{mono-pk}
For $r\geq 5$ and $n\geq 14$,
%
\begin{equation}
p_{r}(n)\geq p_{r}(n-1).
\end{equation}
\end{thm}
It should be noted that  {Theorem~\ref{lem-kernel}} gives more results on
$p_{r}(n)-p_{r}(n-1)$.

The proof of  {Theorem~\ref{main-thm-n}} for $m=2$ and $m=3$ also requires
the monotonicity property of another partition function $pp_{r}(n)$. Let
$pp_{r}(n)$ denote the number of pairs $(\alpha ,\beta )$, where
$\alpha $ is a partition counted by $p_{r}(i)$ and $\beta $ is a partition
counted by $p_{r+1}(n-i)$ for $0\leq i\leq n$. From the definition of
$pp_{r}(n)$ and  {\eqref{gf-pkn}}, it is easy to see that the generating function
of $pp_{r}(n)$ is
%
\begin{equation}
\label{gf-ppkn}
\sum _{n=0}^{\infty }pp_{r}(n)q^{n}=
\frac{1}{(q^{2};q)_{r-1}(q^{2};q)_{r}}.
\end{equation}
We show that $pp_{r}(n)$ has the following monotonicity property.

\begin{thm}%
\label{mono-ppk}
For $r\geq 3$ and $n\geq 8$,
%
\begin{equation}
pp_{r}(n)\geq pp_{r}(n-1).
\end{equation}
\end{thm}

It should be noted that more results on $pp_{r}(n)-pp_{r}(n-1)$ are stated
in  {Theorem~\ref{lem-kernel-2}}.

This paper is organized as follows. In Section~2, we give a brief outline
of the proofs of  {Theorem~\ref{main-un-2-n}} and  {Theorem~\ref{main-thm-n}}. In Section~3, we show the monotonicity property of
$p_{r}(n)$. Section~4 is devoted to the proof of the monotonicity property
of $pp_{r}(n)$. In Section~5, we give a proof of  {Theorem~\ref{main-un-2-n}} by means of {Theorem~\ref{lem-kernel}} and  {Corollary~\ref{corollary}}. Sections~6$\sim $9 are devoted to the proof of  {Theorem~\ref{main-thm-n}}. More specifically, we establish three expressions for
the generating function of $M(m-1,n)-M(m,n)$ in Section~6. Section~7 is
devoted to the proof of  {Theorem~\ref{main-thm-n}} when $m=2$ in light of
 {Theorem~\ref{lem-kernel-2}}. In Section~8, we show  {Theorem~\ref{main-thm-n}} holds when $m\geq 3$ in light of  {Theorem~\ref{lem-kernel}},  {Corollary~\ref{corollary}} and {Theorem~\ref{lem-kernel-2}}. In Section~9, we finish the proof of  {Theorem~\ref{main-thm-n}} by showing that $M(0,n)\geq M(1,n)$ for $n\geq 44$. In
Section~10, we first prove  {Theorem~\ref{lem-indu}}, and then confirm {Conjecture~\ref{conj-opst3}} in light of  {Theorem~\ref{main-thm-n}} and  {Theorem~\ref{lem-indu}}. Finally, two conjectures on the log-concavity
of $M(m,n)$ and $N(m,n)$ are posed in Section~11.

 \section{The outline of the proofs of  {Theorems~\ref{main-un-2-n} and \ref{main-thm-n}}}
\label{sec2}

The proofs of  {Theorem~\ref{main-un-2-n}} and  {Theorem~\ref{main-thm-n}} both
rely on the generating function of $M(m,n)$ established by Garvan
\cite{Garvan-1988-1}:

\begin{thm}{\rm(Garvan).}%
\label{thm-lem-gf-c-1}%
\ For $m\geq 0$,

\begin{equation}
\label{lem-gf-c-1}
\sum _{n= 0}^{\infty }M(m,n)q^{n}=\frac{(1-q)q^{m}}{(q;q)_{m}}+\sum _{k=1}^{\infty }\frac{q^{k(k+m)+2k+m}}{(q;q)_{k}(q^{2};q)_{k+m-1}}.
\end{equation}

\end{thm}

By  {Theorem~\ref{thm-lem-gf-c-1}}, it is easy to see that

\begin{eqnarray}
&&\sum _{n=0}^{\infty }\left (M(m,n)-M(m,n-1)\right ) q^{n}
\nonumber
\\[3pt]
&=& \frac{(1-q)^{2}q^{m}}{(q;q)_{m}}+\frac{q^{2m+3}}{(q^{2};q)_{m}}+
\sum _{k=2}^{\infty
}\frac{q^{k(k+m)+2k+m}}{(q^{2};q)_{k-1}(q^{2};q)_{k+m-1}}.
\label{thm-2-equ-n}
\end{eqnarray}

To prove  {Theorem~\ref{main-un-2-n}}, it suffices to show that the coefficients
of $q^{n}$ in  {\eqref{thm-2-equ-n}} are nonnegative when $n\geq 14$ and
$0\leq m\leq n-2$. By the definition  {\eqref{gf-pkn}} of $p_{r}(n)$ and  {\eqref{thm-2-equ-n}}, we see that for $m\geq 2$ and $n\geq m+1$,

\begin{equation}
\label{equ-m-m-n-m-n-1-2m-3}
M(m,n)-M(m,n-1)\geq p_{m}(n-m)-p_{m}(n-m-1)+p_{m+1}(n-2m-3).
\end{equation}

With the aid of  {Theorem~\ref{lem-kernel}} and  {Corollary~\ref{corollary}}, we will show that the right-hand side of  {\eqref{equ-m-m-n-m-n-1-2m-3}} is nonnegative when $n\geq 14$ and
$0\leq m\leq n-2$, which leads to  {Theorem~\ref{main-un-2-n}}.

Similarly, by  {Theorem~\ref{thm-lem-gf-c-1}}, we see that for
$m\geq 1$,

\begin{eqnarray}
&&\sum _{n= 0}^{\infty }\left ( M(m-1,n)-M(m,n)\right ) q^{n}
\nonumber
\\[3pt]
&=&\sum _{k=1}^{\infty
}\frac{q^{k(k+m-1)+2k+m-1}}{(q;q)_{k}(q^{2};q)_{k+m-2}}-\sum _{k=1}^{\infty }\frac{q^{k(k+m)+2k+m}}{(q;q)_{k}(q^{2};q)_{k+m-1}}
\nonumber
\\[3pt]
&&+\frac{(1-q)q^{m-1}}{(q;q)_{m-1}}- \frac{q^{m}}{(q^{2};q)_{m-1}}.
\label{lem-gf-r-g}
\end{eqnarray}

In order to prove  {Theorem~\ref{main-thm-n}}, we aim to show that the coefficients
of $q^{n}$ in  {\eqref{lem-gf-r-g}} are nonnegative when $n\geq 44$ and
$1\leq m \leq n-1$. It turns out that this will be more difficult and it
is required to transform  {\eqref{lem-gf-r-g}} into several summations which
have nonnegative power series coefficients. To this end, we first split
the first summation in  {\eqref{lem-gf-r-g}} into five summations as stated
in  {Lemma~\ref{lem-m-ge-1-gen}}, and then split the second summation in  {\eqref{lem-gf-r-g}} into five summations as stated in  {Lemma~\ref{lem-m-ge-2-gen}}. Based on  {Lemma~\ref{lem-m-ge-1-gen}} and {Lemma~\ref{lem-m-ge-2-gen}}, we could derive from  {\eqref{lem-gf-r-g}} a new expression
of the generating function of $M(m-1,n)-M(m,n)$ stated in  {Theorem~\ref{lem-4-1}}. Moreover, when $m\ge 2$, it can be shown that some summations in  {Theorem~\ref{lem-4-1}} have nonnegative power series coefficients,  see  {Theorem~\ref{lem-4-2}} and  {Theorem~\ref{lem-4}}.
The proof of  {Theorem~\ref{main-thm-n}} consists of three parts: (1) $m=2$ (Section~\ref{sec7}), (2) $m\ge 3$ (Section~\ref{sec8}), (3) $m=1$ (Section~\ref{sec9}).

 When $m=2$, based on  {Theorem~\ref{lem-4-2}}, we will show that for
$n\geq 15$,
$$M(1,n)-M(2,n)\geq T_{2}(n),$$
where $ T_{2}(n)$ is defined
as:

\begin{equation}
\sum _{n=0}^{\infty }T_{2}(n)q^{n}:=\sum _{r=1}^{\infty
}\frac{q^{r^{2}+7r+7}(1-q)}{(q^{2};q)_{r}(q^{2};q)_{r+1}}.
\end{equation}

By the definition  {\eqref{gf-ppkn}} of $pp_{r}(n)$, we see that

\begin{eqnarray*}
\sum _{n=0}^{\infty }T_{2}(n)q^{n}&=&\sum _{r=1}^{\infty }q^{r^{2}+7r+7}
\sum _{n=0}^{\infty }(pp_{r+1}(n)-pp_{r+1}(n-1))q^{n}.
\end{eqnarray*}

In light of  {Theorem~\ref{lem-kernel-2}}, we will show that
$M(1,n)-M(2,n)\geq 0$ for $n\geq 44$.

When $m\geq 3$. In light of  {Theorem~\ref{lem-4}}, we will show that for
$n\ge 0$,

\begin{equation*}
M(m-1,n)-M(m,n)\ge U_{m}(n),
\end{equation*}

where

\begin{equation*}
\sum _{n=0}^{\infty }U_{m}(n)q^{n}= -q^{2m}+q^{2m+1}+q^{3m+1}+q^{m-1}
\frac{1-q}{(q^{2};q)_{m-2}}+\frac{q^{4m+8}}{(q^{2};q)_{m}}.
\end{equation*}

By the definition  {\eqref{gf-pkn}} of $p_{r}(n)$, we see that
$U_{m}(n)$ can be expressed in terms of $p_{m}(n)$ and
$p_{m}(n)-p_{m}(n-1)$. With the aid of  {Theorem~\ref{lem-kernel}} and {Corollary~\ref{corollary}}, we will show that $U_{m}(n)\ge 0$ for $m\ge 3$ and
$n\ge 44$, which implies $M(m-1,n)\geq M(m,n)$ for $m\ge 3$ and
$n\ge 44$.

The proof of  {Theorem~\ref{main-thm-n}} when $m=1$ is the most complicated.
It is required to do more operations on  {\eqref{lem-gf-r-g}} when
$m=1$. Besides  {Lemma~\ref{lem-m-ge-1-gen}} and  {Lemma~\ref{lem-m-ge-2-gen}}, we also need three more lemmas ( {Lemmas~\ref{main-m-0-po-2-lem}, \ref{main-m-0-po-1-lem}, \ref{lem-0-last}}). Based
on these five lemmas, we succeed to transform  {\eqref{lem-gf-r-g}} when
$m=1$ into several summations stated in  {Theorem~\ref{lem-0-1}} which have
nonnegative power series coefficients. Using  {Theorem~\ref{lem-0-1}}, we
will deduce that $M(0,n)-M(1,n)\ge T_{1}(n)$ for $n\ge 10$, where
$T_{1}(n)$ is defined as in  {\eqref{equ-t1n-q11-q17}}.
With the aid of a result of  Chan and Mao
\cite{Chan-Mao-2014}, namely  {Lemma~\ref{le-chanmao-10}}, and the exact formula of $p_{4}(n)$ stated in {Lemma~\ref{lem-kernel-f}}, we deduce that $T_{1}(n)\ge 0$ for $n\ge 106$. This
leads to $M(0,n)-M(1,n)\ge 0$ for $n\ge 106$. Moreover, it can be checked
that $M(0,n)-M(1,n)\ge 0$ when $44\leq n\leq 105$. Therefore, we show
that  {Theorem~\ref{main-thm-n}} holds when $m=1$.

\section{The monotonicity property of $p_{r}(n)$}
\label{sec3}

In this section, we aim to investigate the monotonicity property of
$p_{r}(n)$. We will prove the following results on
$p_{r}(n)-p_{r}(n-1)$ for $r\geq 2$, which leads to  {Theorem~\ref{mono-pk}} immediately. The results for some
special cases of $p_{r}(n)-p_{r}(n-1)$ will be used in the proofs
of  {Theorem~\ref{main-un-2-n}} and  {Theorem~\ref{main-thm-n}}.

\begin{thm}%
\label{lem-kernel}%
For $r\geq 2$, define
%
\begin{equation}
\label{equ-def-dkn}
d_{r}(n)=p_{r}(n)-p_{r}(n-1).
\end{equation}
Then
\begin{itemize}
\item[\textup{(1)}] $d_{r}(0)=1$ and $d_{r}(1)=-1$ for all $r\ge 2$.
\item[\textup{(2)}] $d_{2}(n)=1$ when $n$ is even and $d_{2}(n)=-1$ when
$n$ is odd.
\item[\textup{(3)}] $d_{3}(n)=1$ when $n\equiv 0,2 \pmod{6}$,
$d_{3}(n)=-1$ when $n \equiv 1 \pmod{6}$ and $d_{3}(n)=0$ when
$n\equiv 3,4,5\pmod{6}$.
\item[\textup{(4)}] $d_{4}(n)> 0$ when $n$ is even,
$d_{4}(n)=-\left \lfloor (n+11)/12\right \rfloor $ when
$n\equiv 1\pmod{2}$ and $n\not\equiv 3 \pmod{12}$ and
$d_{4}(n)=-\left \lfloor n/12\right \rfloor $ when
$n\equiv 3 \pmod{12}$.
\item[\textup{(5)}] $d_{5}(n)\geq 0$ for $n\geq 2$. Moreover,
$d_{5}(n)\geq 1$ for $n\geq 14$.
\item[\textup{(6)}] $d_{6}(n)\geq 0$ for $n\geq 0$ except for
$d_{6}(1)=d_{6}(7)=d_{6}(13)=-1$.
\item[\textup{(7)}] When $r\geq 7$, $d_{r}(n)\geq 0$ for $n\geq 2$. Moreover,
$d_{r}(r+2)\geq 1$ and $d_{r}(2r+7)\geq 1$.
\end{itemize}
\end{thm}

To prove  {Theorem~\ref{lem-kernel}}, we first establish exact formulas for
$p_{r}(n)$ when $2\leq r\leq 4$ in  {Lemma~\ref{lem-kernel-f}}. We then establish
three expressions for the generating function of $p_{r}(n)$ in  {Lemma~\ref{lem-kernel-f2}}. We proceed to show  {Lemma~\ref{lemma-s-1}}, which plays
a crucial role in the proof of  {Theorem~\ref{lem-kernel}}. Finally, we give
a proof of  {Theorem~\ref{lem-kernel}} based on  {Lemma~\ref{lem-kernel-f}},
 {Lemma~\ref{lem-kernel-f2}} and  {Lemma~\ref{lemma-s-1}}.

\begin{lem}
\label{lem-kernel-f}
When $2\leq r\leq 4$, we have the following explicit formulas for
$p_{r}(n)$\textup{:}
\begin{itemize}
\item[\textup{(1)}]
%
\begin{equation}
\label{equ-p2n}
p_{2}(n)=
\begin{cases}
1,&\text{if $n$ is even};
\\
0,&\text{if $n$ is odd.}
\end{cases}
\end{equation}
\item[\textup{(2)}]
%
\begin{align}
p_{3}(n)&=
\begin{cases}
\left \lfloor \frac{n}{6}\right \rfloor +1, & \text{if }\ n\not\equiv 1 \pmod{6};
\\[3pt]
\left \lfloor \frac{n}{6}\right \rfloor , &\text{if }\ n\equiv 1
\pmod{6}.
\end{cases}
\label{lem-2-claim}
\end{align}
\item[\textup{(3)}]
%
\begin{align}
p_{4}(n)&=
\begin{cases}
3a^{2}+3a+1,&\text{if }n=12a \text{\ or\ } n=12a+3;
\\[3pt]
3a^{2}+4a+1,&\text{if }n=12a+2 \text{\ or\ } n=12a+5;
\\[3pt]
3a^{2}+5a+2,&\text{if }n=12a+4 \text{\ or\ } n=12a+7;
\\[3pt]
3a^{2}+6a+3,&\text{if }n=12a+6 \text{\ or\ } n=12a+9;
\\[3pt]
3a^{2}+7a+4,&\text{if }n=12a+8 \text{\ or\ } n=12a+11;
\\[3pt]
3a^{2}+8a+5,&\text{if }n=12a+10 \text{\ or\ } n=12a+13.
\end{cases}
\label{lem-3-claim}
\end{align}
\end{itemize}
\end{lem}

\begin{proof} Let $p(n,r)$ denote the number of partitions of $n$ with at
most $r$ parts. It is well known that
%
\begin{equation}
\label{equ-gen-pnkqn}
\sum _{n=0}^{\infty }p(n,r)q^{n}=\frac{1}{(q;q)_{r}}.
\end{equation}

From  {\eqref{gf-pkn}} and  {\eqref{equ-gen-pnkqn}}, we see that
%
\begin{equation}
\sum _{n=0}^{\infty }p_{r}(n)q^{n}=\frac{1-q}{(q;q)_{r}}=\sum _{n=0}^{\infty }p(n,r)q^{n}-\sum _{n=1}^{\infty }p(n-1,r)q^{n}.
\end{equation}
Therefore, we find that for $n\ge 1$,
%
\begin{equation}
\label{equ-pn2-left-final}
p_{r}(n)=p(n,r)-p(n-1,r).
\end{equation}

When $r=2$, Andrews \cite{Andrews-2003} showed that
%
\begin{equation}
\label{equ-pn2-left}
p(n,2)=\left \lfloor \frac{n+2}{2}\right \rfloor .
\end{equation}
DeMorgan \cite{DeMorgan-1843} found the following formula for
$p(n,3)$ as given below,
%
\begin{equation}
\label{equ-pn3-left}
p(n,3)=\left \{  \frac{(n+3)^{2}}{12}\right \}  .
\end{equation}
Gl\"{o}sel \cite{Glosel-1896} gave the following formula for
$p(n,4)$,
%
\begin{equation}
\label{equ-pn4-left}
p(n,4)=\left \{  \left \lfloor \frac{(n+4)}{2}\right \rfloor ^{2}
\left (3\left \lfloor \frac{n+9}{2}\right \rfloor -\left \lfloor \frac{n+10}{2} \right \rfloor \right )\frac{1}{36}\right \}  ,
\end{equation}
where $\lfloor x\rfloor $ is the greatest integer $\le x$, and
$\{x\}$ is the nearest integer to $x$. Substituting  {\eqref{equ-pn2-left}}, {\eqref{equ-pn3-left}} and {\eqref{equ-pn4-left}} into  {\eqref{equ-pn2-left-final}}, and after some calculations, we see that {\eqref{equ-p2n}},  {\eqref{lem-2-claim}} and  {\eqref{lem-3-claim}} hold.
\end{proof}

By  {Lemma~\ref{lem-kernel-f}}, we obtain the following corollary, which
is useful in the proof of  {Theorem~\ref{main-un-2-n}}.

\begin{core}%
\label{corollary}
For $r\geq 3$ and $n\geq 2$, we have $p_{r}(n)\geq 1$. Moreover,
$p_{r}(n)\geq \left \lfloor \frac{n}{6}\right \rfloor $.
\end{core}
 \begin{proof} By the definition of $p_{r}(n)$, it is clear to see that for
any $i\geq 2$, each partition counted by $p_{i}(n)$ is also counted by
$p_{i+1}(n)$. So
%
\begin{equation}
\label{coro-temp}
p_{r}(n)\geq p_{r-1}(n)\geq \cdots \geq p_{3}(n).
\end{equation}
Furthermore,  it is clear from  {Lemma~\ref{lem-kernel-f}} that $p_{3}(n)\geq 1$ for
$n\geq 2$. Moreover,
$p_{3}(n)\geq \left \lfloor \frac{n}{6}\right \rfloor $. This yields the
corollary. \end{proof}

The following lemma gives three expressions for the generating function
of $p_{r}(n)$. To be specific, the expression  {\eqref{lem-snk-gf-3}} will
be used in the proof of  {Theorem~\ref{lem-4}} and the expressions  {\eqref{lem-snk-gf-1}} and {\eqref{lem-snk-gf-2}} will be used in the proof
of  {Theorem~\ref{lem-kernel}}.

\begin{lem}
\label{lem-kernel-f2}
For $r\ge 2$. we have
%
\begin{align}
\sum _{n=0}^{\infty }p_{r}(n)q^{n}&=1+\sum _{j=2}^{r}
\frac{q^{j}}{(q^{j};q)_{r-j+1}}
\label{lem-snk-gf-3}
\\[3pt]
&=1-q+\frac{q}{(q^{2};q)_{r-2}} +\sum _{j=1}^{r}
\frac{q^{2j}}{(q^{2};q)_{j-1}}
\label{lem-snk-gf-1}
\\[3pt]
&=q^{r}+\frac{1}{(q^{2};q)_{r-2}}+\frac{q^{2r}}{(q^{2};q)_{r-1}}+
\sum _{j=2}^{r-1} \frac{q^{r+j}}{(q^{2};q)_{j-1}}.
\label{lem-snk-gf-2}
\end{align}
\end{lem}

\begin{proof} We first verify  {\eqref{lem-snk-gf-3}}. For
$2\leq j\leq r$, let $p_{r,j}(n)$ denote the number of partitions of
$n$ such that each part is not exceeding $r$ and the smallest part is equal
to $j$. Clearly, for $n\geq 1$,
%
\begin{equation}
\label{equ-pkjn-pkn}
\sum _{j=2}^{r} p_{r,j}(n)=p_{r}(n).
\end{equation}
On the other hand, it is easy to see that
%
\begin{equation}
\label{equ-pkn-cls-sml}
\sum _{n=1}^{\infty }p_{r,j}(n)q^{n}=\frac{q^{j}}{(q^{j};q)_{r-j+1}}.
\end{equation}
Combining  {\eqref{equ-pkjn-pkn}} and  {\eqref{equ-pkn-cls-sml}}, we obtain
\begin{equation*}
\sum _{n=0}^{\infty }p_{r}(n)q^{n}=1+\sum _{j=2}^{r}\sum _{n= 1}^{\infty }p_{r,j}(n)q^{n}=1+\sum _{j=2}^{r}
\frac{q^{j}}{(q^{j};q)_{r-j+1}},
\end{equation*}
which is  {\eqref{lem-snk-gf-3}}.

We proceed to derive  {\eqref{lem-snk-gf-1}}. To this end, we need to divide
the set of partitions counted by $p_{r}(n)$ into two disjoint subsets based
on the difference of the largest part of the partition and the second largest
part. Let $s_{r}(n)$ denote the number of partitions
$\lambda =(\lambda _{1},\lambda _{2},\ldots ,\lambda _{\ell})$ counted by
$p_{r}(n)$ such that $\lambda _{1}-\lambda _{2}\geq 1 $ and
$q_{r}(n)$ denote the number of partitions
$\lambda =(\lambda _{1},\lambda _{2},\ldots ,\lambda _{\ell})$ counted by
$p_{r}(n)$ such that $\lambda _{1}-\lambda _{2}=0$. Here we use the convention
that $\lambda _{i}=0$ for $i>\ell$. Obviously,
\begin{equation*}
p_{r}(n)=s_{r}(n)+q_{r}(n).
\end{equation*}
Hence the generating function of $p_{r}(n)$ is equal to the sum of the
generating functions of $s_{r}(n)$ and $q_{r}(n)$.

We first consider the generating function of $s_{r}(n)$. Let
$\lambda =(\lambda _{1}, \lambda _{2},\ldots , \lambda _{\ell})$ be a partition
counted by $s_{r}(n)$. If $\lambda \neq (2)$ note that
$\lambda _{1}>\lambda _{2}$. Then we can define
$\mu =(\lambda _{1}-1,\lambda _{2},\ldots , \lambda _{\ell})$ which clearly
is a partition counted by $p_{r-1}(n-1)$. Hence, by  {\eqref{gf-pkn}}, we
obtain the following generating function of $s_{r}(n)$:
%
\begin{equation}
\label{lem-snk-gf-1-case1}
\sum _{n=0}^{\infty }s_{r}(n)q^{n} =-q+q^{2}+\frac{q}{(q^{2};q)_{r-2}}.
\end{equation}
To establish the generating function of $q_{r}(n)$, we will classify the
set of partitions counted by $q_{r}(n)$ based on the size of the largest
part. Given $2\leq j\leq r$, let $q_{r,j}(n)$ denote the number of
partitions counted by $q_{r}(n)$ with the largest part $j$. By definition,
we see that the generating function of $q_{r,j}(n)$ is equal to
%
\begin{align}
\label{lem-snk-gf-1-tem214}
\sum _{n=0}^{\infty }q_{r,j}(n)q^{n} &=\frac{1}{1-q^{2}}
\frac{1}{1-q^{3}}\cdots \frac{1}{1-q^{j-1}}\frac{q^{2j}}{1-q^{j}}=
\frac{q^{2j}}{(q^{2};q)_{j-1}}.
\end{align}
Notice that the empty partition of $0$ is counted by $q_{r}(n)$, so
%
\begin{equation}
\label{lem-snk-gf-1-case2}
\sum _{n=0}^{\infty }q_{r}(n)q^{n}=1+\sum _{j=2}^{r}\sum _{n= 0}^{\infty }q_{r,j}(n)q^{n} =1+\sum _{j=2}^{r}
\frac{q^{2j}}{(q^{2};q)_{j-1}}.
\end{equation}
Combining  {\eqref{lem-snk-gf-1-case1}} and  {\eqref{lem-snk-gf-1-case2}}, we
obtain  {\eqref{lem-snk-gf-1}}.

We finish the proof of  {Lemma~\ref{lem-kernel-f}} by showing  {\eqref{lem-snk-gf-2}} holds. We first divide the set of partitions counted
by $p_{r}(n)$ into three disjoint sets. Let $g_{r}(n)$ denote the number
of partitions
$\lambda =(\lambda _{1},\lambda _{2},\ldots ,\lambda _{\ell})$ counted by
$p_{r}(n)$ with $r=\lambda _{1}=\lambda _{2}$ and $h_{r}(n)$ denote the
number of partitions
$\lambda =(\lambda _{1},\lambda _{2},\ldots ,\lambda _{\ell})$ counted by
$p_{r}(n)$ with $r=\lambda _{1}>\lambda _{2}$. Note that the number of
partitions
$\lambda =(\lambda _{1},\lambda _{2},\ldots ,\lambda _{\ell})$ counted by
$p_{r}(n)$ with $r>\lambda _{1}$ is equal to $p_{r-1}(n)$, hence
\begin{equation*}
p_{r}(n)=p_{r-1}(n)+g_{r}(n)+h_{r}(n).
\end{equation*}
By  {\eqref{gf-pkn}}, we see that the generating function of
$p_{r-1}(n)$ is
%
\begin{equation}
\label{lem-snk-gf-2-case1}
\sum _{n= 0}^{\infty }p_{r-1}(n)q^{n}=\frac{1}{(q^{2};q)_{r-2}}.
\end{equation}
Notice that $g_{r}(n)$ coincides with $q_{r,r}(n)$. Hence by  {\eqref{lem-snk-gf-1-tem214}}, we see that
%
\begin{equation}
\label{lem-snk-gf-2-case2}
\sum _{n= 0}^{\infty }g_{r}(n)q^{n}=\frac{q^{2r}}{(q^{2};q)_{r-1}}.
\end{equation}
To obtain the generating function of $h_{r}(n)$, we define
$h_{r,j}(n)$ as the number of partitions
$\lambda =(\lambda _{1},\lambda _{2},\ldots ,\lambda _{\ell})$ counted by
$h_{r}(n)$ with the second largest part $\lambda _{2}=j$. If $j=0$, then
$\lambda =(r)$. Otherwise, for $2\leq j\leq r-1$, the generating function
of $h_{r,j}(n)$ is equal to
%
\begin{align}
\label{lem-snk-gf-1-tem}
\sum _{n= 0}^{\infty }h_{r,j}(n)q^{n} &=\frac{1}{1-q^{2}}
\frac{1}{1-q^{3}}\cdots \frac{1}{1-q^{j-1}}\frac{q^{j}}{1-q^{j}}q^{r}=
\frac{q^{r+j}}{(q^{2};q)_{j-1}}.
\end{align}
Hence, we obtain the following generating function of $h_{r}(n)$
%
\begin{equation}
\label{lem-snk-gf-2-case3}
\sum _{n= 0}^{\infty }h_{r}(n)q^{n}=q^{r}+\sum _{j=2}^{r-1}\sum _{n= 0}^{\infty }h_{r,j}(n)q^{n}=q^{r}+\sum _{j=2}^{r-1}
\frac{q^{r+j}}{(q^{2};q)_{j-1}}.
\end{equation}
Combining  {\eqref{lem-snk-gf-2-case1}},  {\eqref{lem-snk-gf-2-case2}} and {\eqref{lem-snk-gf-2-case3}}, we obtain  {\eqref{lem-snk-gf-2}}. This completes
the proof. \end{proof}

Before proceeding to prove the monotonicity property of $p_{r}(n)$, let us first show {Lemma~\ref{lemma-s-1}}, which plays a crucial role
in the proof of  {Theorem~\ref{lem-kernel}}. It turns out that  the following lemma due to Chan
and Mao \cite{Chan-Mao-2014} is needed in the proof of
 {Lemma~\ref{lemma-s-1}}.

\begin{lem}{\rm(Chan and Mao).}%
\label{le-chanmao-10}
\begin{equation*}
\frac{1-q^{m}}{(1-q^{2})(1-q^{3})}
\end{equation*}
has nonnegative power series coefficients for any integer $m\ge 2$.
\end{lem}

\begin{lem}%
\label{lemma-s-1}
For $r\geq 4$, let
%
\begin{equation}
\label{lemma-s-1-e}
\sum _{n=0}^{\infty }t_{r}(n)q^{n}:=\sum _{j=2}^{r}
\frac{q^{2j}(1-q^{r-j+2})}{(q^{2};q)_{j-1}}.
\end{equation}
Then $t_{r}(n)\geq 0$ for $n\geq 0$. Moreover, when $r\neq 5$, we have
$t_{r}(n)\geq 1$ for $n\geq 14$.
\end{lem}

\begin{proof} Define
\begin{equation*}
\sum _{n= 0}^{\infty }t^{(j)}_{r}(n)q^{n}=
\frac{q^{2j}(1-q^{r-j+2})}{(q^{2};q)_{j-1}},
\end{equation*}
obviously,
%
\begin{equation}
\label{lemma-s-1-e-sum}
t_{r}(n)=\sum _{j=2}^{r}t^{(j)}_{r}(n).
\end{equation}
By  {Lemma~\ref{le-chanmao-10}}, we see that when $r\geq 4$ and
$3\leq j\leq r$,
\begin{equation*}
\sum _{n=0}^{\infty }t^{(j)}_{r}(n)q^{n}=
\frac{q^{2j}(1-q^{r-j+2})}{(q^{2};q)_{j-1}}=
\frac{1-q^{r-j+2}}{(1-q^{2})(1-q^{3})}\cdot
\frac{q^{2j}}{(q^{4};q)_{j-3}}
\end{equation*}
has nonnegative power series coefficients. It gives that when
$r\geq 4$ and $3\leq j\leq r$,
%
\begin{equation}
\label{lemma-s-1-e-jgeq3}
t^{(j)}_{r}(n)\geq 0 \quad \text{for} \quad n\geq 0.
\end{equation}
We next show that $t^{(2)}_{r}(n)+t^{(3)}_{r}(n)\geq 0$ when
$r\neq 5$. First, it is easy to see that
%
\begin{align}
\sum _{n=0}^{\infty }t^{(2)}_{r}(n)q^{n}+\sum _{n=0}^{\infty }t^{(3)}_{r}(n)q^{n}&=
\frac{q^{4}(1-q^{r})}{1-q^{2}} +
\frac{q^{6}(1-q^{r-1})}{(1-q^{2})(1-q^{3})}
\nonumber
\\[3pt]
&=
\frac{q^{4}-q^{7}+q^{6}-q^{r+2}+q^{r+2}-q^{r+4}-q^{r+5}+q^{r+7}}{(1-q^{2})(1-q^{3})}
\nonumber
\\[3pt]
&=\frac{q^{4}}{1-q^{2}}+\frac{q^{6}-q^{r+2}}{(1-q^{2})(1-q^{3})}+q^{r+2}.
\label{lemma-s-1-e-j23-temp}
\end{align}
By using  {Lemma~\ref{le-chanmao-10}} again, we find that when
$r\geq 4$ and $r\neq 5$,
\begin{equation*}
\frac{q^{6}-q^{r+2}}{(1-q^{2})(1-q^{3})}=
\frac{q^{6}(1-q^{r-4})}{(1-q^{2})(1-q^{3})}
\end{equation*}
has nonnegative power series coefficients. Hence, from  {\eqref{lemma-s-1-e-j23-temp}}, we see that when $r\neq 5$,
%
\begin{equation}
\label{lemma-s-1-e-j23}
t^{(2)}_{r}(n)+t^{(3)}_{r}(n)\geq 0 \quad \text{for} \quad n\geq 0.
\end{equation}
Thus, we derive from  {\eqref{lemma-s-1-e-sum}} and  {\eqref{lemma-s-1-e-jgeq3}} that $t_{r}(n)\geq 0$ when $r\neq 5$.

We next show that when $r\neq 5$, $t_{r}(n)\geq 1$ for $n\geq 14$. By the
generating function  {\eqref{gf-pkn}} of $p_{r}(n)$, we see that
%
\begin{equation}
\label{equ-lemma-s-1-2}
\sum _{n=0}^{\infty }t^{(4)}_{r}(n)q^{n}=
\frac{q^{8}(1-q^{r-2})}{(1-q^{2})(1-q^{3})(1-q^{4})}=\sum _{n=8}^{\infty }(p_{4}(n-8)-p_{4}(n-r-6))q^{n}.
\end{equation}
From  {Lemma~\ref{lem-kernel-f}} (3), it is easy to check that for
$n\geq 14$,
\begin{equation*}
p_{4}(n-8)>p_{4}(n-10),
\end{equation*}
and
\begin{equation*}
p_{4}(n-8)>p_{4}(n-13).
\end{equation*}
Hence when $r\ge 4$ even,
\begin{equation*}
p_{4}(n-8)>p_{4}(n-10)\geq \cdots \geq p_{4}(n-r-6),
\end{equation*}
and when $r\geq 7$ odd,
\begin{equation*}
p_{4}(n-8)>p_{4}(n-13)\geq p_{4}(n-15)\geq \cdots \geq p_{4}(n-r-6).
\end{equation*}
In either case, we see that $p_{4}(n-8)>p_{4}(n-r-6)$ for $n\geq 14$. It
yields that when $r\neq 5$,
%
\begin{equation}
\label{lemma-s-1-e-posit}
t^{(4)}_{r}(n)\geq 1 \quad \text{for} \quad n\geq 14.
\end{equation}
Combining  {\eqref{lemma-s-1-e-jgeq3}},  {\eqref{lemma-s-1-e-j23}} and {\eqref{lemma-s-1-e-posit}}, it follows from  {\eqref{lemma-s-1-e-sum}} that
when $r\neq 5$, $t_{r}(n)\geq 1$ for $n\geq 14$. Thus, we complete the
proof of  {Lemma~\ref{lemma-s-1}} when $r\neq 5$.

It remains to show that  {Lemma~\ref{lemma-s-1}} holds when $r=5$. From {\eqref{lemma-s-1-e-jgeq3}}, we see that for $3\leq j\leq 5$,
%
\begin{equation}
\label{lemma-s-1-e-jgeq34}
t^{(j)}_{5}(n)\geq 0.
\end{equation}
Note that
\begin{equation*}
\sum _{n=0}^{\infty }t^{(2)}_{5}(n)q^{n}=
\frac{q^{4}(1-q^{5})}{(1-q^{2})}= \frac{q^{4}}{1-q^{2}}-
\frac{q^{9}}{1-q^{2}}=\sum _{n= 2}^{\infty }q^{2n}-\sum _{n= 4}^{\infty }q^{2n+1},
\end{equation*}
so we derive that $t^{(2)}_{5}(n)=1$ when $n$ is even and $n\geq 4$ and
$ t^{(2)}_{5}(n)=-1$ when $n$ is odd and $n\geq 9$. Since
\begin{equation*}
\sum _{n=0}^{\infty }t^{(5)}_{5}(n)q^{n}= \frac{q^{10}}{(q^{3};q)_{3}}=
\sum _{i=0}^{\infty }\sum _{j=0}^{\infty }\sum _{h=0}^{\infty }q^{3i+4j+5h+10},
\end{equation*}
and it is easy to check that for $n\geq 13$, {there exists nonnegative} integers
$i,j,h$ such that $3i+4j+5h+10=n$. Hence $t^{(5)}_{5}(n)\geq 1$. Thus,
from the above analysis, we derive that
\begin{equation*}
t_{5}(n)\geq t^{(2)}_{5}(n)+t^{(5)}_{5}(n)\geq 0 \quad \text{for}
\quad n\geq 13.
\end{equation*}
It is trivial to check that $t_5(n)\ge 0$ for $0\le n\le 12$. Hence  {Lemma~\ref{lemma-s-1}} is also valid when $r=5$. Thus, we complete
the proof of the lemma. \end{proof}

We are now in a position to give a proof of  {Theorem~\ref{lem-kernel}}.

\begin{proof}[Proof of  {Theorem~\ref{lem-kernel}}] From  {\eqref{equ-def-dkn}} and
 {Lemma~\ref{lem-kernel-f}}, it is easy to check that  {Theorem~\ref{lem-kernel}} holds when $r=2,3$, or $4$.

We now consider the case $r\geq 5$. By  {\eqref{lem-snk-gf-1}}, we see that
\begin{equation*}
\sum _{n=0}^{\infty }p_{r}(n)q^{n}=1-q+\frac{q}{(q^{2};q)_{r-2}} +\sum _{j=1}^{r}
\frac{q^{2j}}{(q^{2};q)_{j-1}},
\end{equation*}
and from  {\eqref{lem-snk-gf-2}},
\begin{equation*}
\sum _{n=1}^{\infty }p_{r}(n-1)q^{n}= q^{r+1}+\frac{q}{(q^{2};q)_{r-2}}+
\frac{q^{2r+1}}{(q^{2};q)_{r-1}} +\sum _{j=2}^{r-1}
\frac{q^{r+j+1}}{(q^{2};q)_{j-1}}.
\end{equation*}
Hence, we have the following generating function
of $d_{r}(n)$:
\begin{align*}
\sum _{n=0}^{\infty }d_{r}(n)q^{n}&=1-q-q^{r+1} -
\frac{q^{2r+1}}{(q^{2};q)_{r-1}}+\sum _{j=1}^{r}
\frac{q^{2j}}{(q^{2};q)_{j-1}} -\sum _{j=2}^{r-1}
\frac{q^{r+j+1}}{(q^{2};q)_{j-1}},
\end{align*}
which can be simplified as
%
\begin{equation}
\label{lem-4-39}
\sum _{n= 0}^{\infty }d_{r}(n)q^{n}=1-q+q^{2}-q^{r+1}+q^{2r}
\frac{1-q}{(q^{2};q)_{r-1}}+\sum _{j=2}^{r-1}
\frac{q^{2j}(1-q^{r-j+1})}{(q^{2};q)_{j-1}} .
\end{equation}
From  {\eqref{gf-pkn}} and  {\eqref{equ-def-dkn}}, it is easy to see that
\begin{equation*}
q^{2r}\frac{1-q}{(q^{2};q)_{r-1}}=\sum _{n=2r}^{\infty }d_{r}(n-2r)q^{n}.
\end{equation*}
Moreover, using the notation of $t_{r}(n)$ as defined in  {\eqref{lemma-s-1-e}}, we see that {\eqref{lem-4-39}} can be expressed as
%
\begin{equation}
\label{thm-6-1}
\sum _{n=0}^{\infty }d_{r}(n)q^{n}=1-q+q^{2}-q^{r+1} +\sum _{n=2}^{\infty }t_{r-1}(n)q^{n} +\sum _{n= 2r}^{\infty }d_{r}(n-2r)q^{n}.
\end{equation}
Hence we obtain the following recurrence relation:
%
\begin{equation}
\label{equ-dkn-tk-1n-ns}
d_{r}(n)=
\begin{cases}
t_{r-1}(n)+1,&\text{if }n=2\text{ or }2r;
\\[3pt]
t_{r-1}(n)-1,&\text{if }n=r+1\text{ or }2r+1;
\\[3pt]
t_{r-1}(n),& \text{if } 3\leq n\leq 2r-1\text{ and } n\neq r+1;
\\[3pt]
t_{r-1}(n)+d_{r}(n-2r), & \text{if } n\geq 2r+2.
\end{cases}
\end{equation}

We next show that  {Theorem~\ref{lem-kernel}} holds when $r\geq 5$. From  {Lemma~\ref{lemma-s-1}}, we see that  $t_{r-1}(n)\geq 0$ for
$r\ge 5$ and $n\geq 2$. Moreover, when $r=5$ or $r\geq 7$, $t_{r-1}(n)\geq 1$ for
$n\geq 14$, which implies that $t_{r-1}(r+1)\geq 1$ and
$t_{r-1}(2r+1)\geq 1$ when $r\geq 13$. By a simple calculation, we find
that $t_{r-1}(r+1)\geq 1$ for $5\leq r\leq 12$ and $r\neq 6$; and
$t_{4}(11)\geq 1$. It follows that $t_{r-1}(r+1)\geq 1$ and
$t_{r-1}(2r+1)\geq 1$ when $r=5$ or $r\geq 7$. Hence, by  {\eqref{equ-dkn-tk-1n-ns}}, we derive that when $r=5$ or $r\geq 7$,
%
\begin{equation}
\label{equ-dkn-tk-1n-nsaa}
d_{r}(n)\geq 0 \quad \text{ for } 2\leq n\leq 2r+1,
\end{equation}
and when $r=6$,
%
\begin{equation}
\label{equ-drn-tk-1n-nsaa}
d_{6}(n)\geq 0 \quad \text{ for } 2\leq n\leq 12\text{ and } n\neq 7.
\end{equation}
We proceed to show that when $r\geq 5$, $d_{r}(n)\geq 0$ for
$n\geq 2r+2$ by using induction on $n$. Assume that there exists a positive
integer $N_{r}\geq 2r+1$ such that when $r=5$ or $r\geq 7$,
$d_{r}(n)\geq 0$ for $2\leq n\leq N_{r}$ and $d_{6}(n)\geq 0$ for
$2\leq n\leq N_{6}$ and $n\neq 7,13$.

We proceed to show that $d_{r}(N_{r}+1)\geq 0$ when $r\geq 5$. By  {\eqref{equ-dkn-tk-1n-ns}} and the fact that $N_{r}+1\geq 2r+2$, we have
%
\begin{equation}
\label{equ-dkn-tk-1n-nstt}
d_{r}(N_{r}+1)=t_{r-1}(N_{r}+1)+d_{r}(N_{r}-2r+1).
\end{equation}
From  {Lemma~\ref{lemma-s-1}}, we see that $t_{r-1}(N_{r}+1)\geq 0$ when
$r\geq 5$ and by the induction hypothesis, we see that when $r=5$ or
$r\geq 7$, $d_{r}(N_{r}+1-2r)\geq 0$ and $d_{6}(N_{6}-11)\geq 0$ for
$N_{6}\neq 18$ or $24$. Hence we derive from  {\eqref{equ-dkn-tk-1n-nstt}} that when $r=5$ or $r\geq 7$,
$d_{r}(N_{r}+1)\geq 0$ and $d_{6}(N_{6}+1)\geq 0$ for $N_{6}\neq 18$ or
$24$. Moreover, it is easy to check that $d_{6}(19)\geq 0$ and
$d_{6}(25)\geq 0$. Thus, we conclude that when $r=5$ or $r\geq 7$,
$d_{r}(n)\geq 0$ for $n\geq 2$ and $d_{6}(n)\geq 0$ for $n\geq 14$.

We finish the proof of  {Theorem~\ref{lem-kernel}} by considering the positivity
of $d_{r}(n)$ when $r=5$ or $r\geq 7$. When $r=5$, and by  {\eqref{equ-dkn-tk-1n-ns}}, we see that when $n\geq 14$,
\begin{equation*}
d_{5}(n)=t_{4}(n)+d_{5}(n-10).
\end{equation*}
Since $t_{4}(n)\geq 1$ for $n\geq 14$ and $d_{5}(n-10)\geq 0$ for
$n\geq 14$, we deduce that $d_{5}(n)\geq 1$ for $n\geq 14$. Thus we complete
the proof of  {Theorem~\ref{lem-kernel}} when $r=5$.

From  {Lemma~\ref{lemma-s-1}}, we see that $t_{r-1}(n)\geq 1$ for
$n\geq 14$ and $r\ge 7$. It follows from  {\eqref{equ-dkn-tk-1n-ns}} that
$d_{r}(r+2)=t_{r-1}(r+2)\geq 1$ when $r\geq 12$. Furthermore, it is easy
to check that $d_{r}(r+2)\geq 1$ when $7\leq r\leq 11$. So
$d_{r}(r+2)\geq 1$ when $r\geq 7$. On the other hand, by  {\eqref{equ-dkn-tk-1n-ns}}, we see that
\begin{equation*}
d_{r}(2r+7)=t_{r-1}(2r+7)+d_{r}(7).
\end{equation*}
Note that $t_{r-1}(2r+7)\geq 1$ when $r\geq 7$ and $d_{r}(7)\geq 0$ when
$r\geq 7$, so we arrive at $d_{r}(2r+7)\geq 1$. Thus we complete the proof
of  {Theorem~\ref{lem-kernel}}.\end{proof}

\section{The monotonicity property of $pp_{r}(n)$}
\label{sec4}

This section is devoted to the monotonicity property of $pp_{r}(n)$ in
terms of  {Theorem~\ref{lem-kernel}}. We will show the following results on
$pp_{r}(n)-pp_{r}(n-1)$, which gives  {Theorem~\ref{mono-ppk}} immediately.
The results for some special cases of
$pp_{r}(n)-pp_{r}(n-1)$ will   be used in the proof of  {Theorem~\ref{main-thm-n}} when $m=2$.

\begin{thm}%
\label{lem-kernel-2}%
For $r\geq 2$, let
\begin{equation*}
f_{r}(n)=pp_{r}(n)-pp_{r}(n-1).
\end{equation*}
Then
\begin{itemize}
\item[\textup{(1)}] $f_{r}(0)=1$ and $f_{r}(1)=-1$.
\item[\textup{(2)}] $f_{2}(n)\geq 0$ if and only if $n$ is even. Moreover,
$f_{2}(n)= -\left \lceil \frac{n}{6}\right \rceil $ when $n$ is odd.
\item[\textup{(3)}] $f_{3}(n)\geq 0$ for $n\geq 2$ and $n\neq 7$. Moreover,
$f_{3}(3)=f_{3}(5)=0$, $f_{3}(7)=-1$ and $f_{3}(n)\geq (n-15)/2$ when
$n$ is odd and $n\geq 17$.
\item[\textup{(4)}] When $r\geq 4$, $f_{r}(n)\geq 0$ for $n\geq 2$. Moreover,
$f_{r}(2r+7)\geq 1$.
\end{itemize}
\end{thm}
 \begin{proof} (1) From the definition of $f_{r}(n)$, and by  {\eqref{gf-ppkn}}, we see that
%
\begin{equation}
\label{lem-ker-2-equ}
\sum _{n=0}^{\infty }f_{r}(n)q^{n}=
\frac{1-q}{(q^{2};q)_{r-1}(q^{2};q)_{r}}.
\end{equation}
Clearly $f_{r}(0)=1$ and $f_{r}(1)=-1$.

(2) When $r=2$, we see that
\begin{eqnarray*}
\sum _{n=0}^{\infty }f_{2}(n)q^{n} =\frac{1-q}{(1-q^{2})(q^{2};q)_{2}}=
\frac{1}{1-q^{2}}\sum _{n= 0}^{\infty }d_{3}(n)q^{n}.
\end{eqnarray*}
From  {Theorem~\ref{lem-kernel}} (2), we find that
%
\begin{equation}
\label{lem-2-kernel-1-equ}
\sum _{n=0}^{\infty }f_{2}(n)q^{n} =\frac{1}{1-q^{2}}\sum _{i=0}^{\infty }(q^{6i}-q^{6i+1}+q^{6i+2}),
\end{equation}
which implies that
\begin{equation*}
\sum _{n=0}^{\infty }f_{2}(2n)q^{2n}=\frac{1}{1-q^{2}}\sum _{i=0}^{\infty }(q^{6i}+q^{6i+2}),
\end{equation*}
and
\begin{equation*}
\sum _{n=0}^{\infty }f_{2}(2n+1)q^{2n+1}=-\frac{1}{1-q^{2}}\sum _{i=0}^{\infty }q^{6i+1}=-\sum _{i=0}^{\infty }\sum _{j=0}^{\infty }q^{6i+2j+1}.
\end{equation*}
Hence $f_{2}(2n)\geq 0$ and
$f_{2}(2n+1)=-\left \lceil \frac{2n+1}{6}\right \rceil $ for
$n\geq 0$. Thus we complete the proof of  {Theorem~\ref{lem-kernel-2}} when
$r=2$.

(3) When $r=3$, we see that
%
\begin{eqnarray}%
\label{fk3}
\sum _{n= 0}^{\infty }f_{3}(n)q^{n} &=&
\frac{1-q}{(q^{2};q)_{2}(q^{2};q)_{3}}
\nonumber
\\[3pt]
&=&\frac{1-q}{(q^{2};q)_{4}}\cdot \frac{1-q^{5}}{(1-q^{2})(1-q^{3})}
\nonumber
\\[3pt]
&=&\left (\sum _{n= 0}^{\infty }d_{5}(n)q^{n}\right )\left (
\frac{1}{1-q^{3}} +\frac{q^{2}}{1-q^{2}}\right ).
\end{eqnarray}

Define
\begin{eqnarray*}%
\sum _{n=0}^{\infty }s(n)q^{n}=\left (1-q+q^{4}+\sum_{{n\geq 15}\atop {n
\text{ odd}}}q^{n}\right )\left (\frac{1}{1-q^{3}} +
\frac{q^{2}}{1-q^{2}}\right ).
\end{eqnarray*}
From  {Theorem~\ref{lem-kernel}} (5), we see that $d_{5}(n)\geq 1$ for
$n\geq 14$ and note that $d_{5}(4)\geq 1$. Hence, we deduce that for
$n\geq 0$,
%
\begin{equation}
\label{f3g}
f_{3}(n)\geq s(n).
\end{equation}
Observe that
%
\begin{eqnarray}%
\label{f3pfg}
\sum _{n=0}^{\infty }s(n)q^{n}=\left (1-q+q^{4}+q^{15}\right )\left (
\frac{1}{1-q^{3}} +\frac{q^{2}}{1-q^{2}}\right )+\sum _{{n\geq 17
}\atop{n \text{ odd}}}q^{n}\left (\frac{1}{1-q^{3}} +
\frac{q^{2}}{1-q^{2}}\right ),\nonumber\\
\end{eqnarray}
and note that
\begin{eqnarray*}
&&(1-q+q^{4}+q^{15}) \left (\frac{1}{1-q^{3}}+\frac{q^{2}}{1-q^{2}}
\right )
\\[3pt]
&=&\frac{1+q^{15}}{1-q^{3}}+\frac{-q+q^{4}}{1-q^{3}}+
\frac{q^{2}(1+q^{4})}{1-q^{2}}+\frac{q^{2}(-q+q^{15})}{1-q^{2}}
\\[3pt]
&=&\frac{1+q^{15}}{1-q^{3}}+\frac{q^{2}+q^{6}}{1-q^{2}}-q-(q^{3}+q^{5}+
\cdots +q^{15}),
\end{eqnarray*}
so we find that $s(n)\geq 0$ for $n\geq 16$. Moreover,
\begin{eqnarray*}
\sum _{{n\geq 17}\atop{n \text{ odd}}}q^{n}\left (\frac{1}{1-q^{3}}+
\frac{q^{2}}{1-q^{2}}\right ) &=&\sum _{n=8}^{\infty }q^{2n+1}\left (1+
\frac{q^{3}}{1-q^{3}}+\frac{q^{2}}{1-q^{2}}\right )
\\[3pt]
&=&\frac{1}{1-q^{2}}\sum _{n=8}^{\infty }q^{2n+1}+\frac{q^{3}}{1-q^{3}}
\sum _{n= 8}^{\infty }q^{2n+1}
\\[3pt]
&=&\frac{q^{17}}{(1-q^{2})^{2}}+\frac{q^{3}}{1-q^{3}}\sum _{n=8}^{\infty }q^{2n+1}
\\[3pt]
&=&\sum _{{n\geq 17}\atop{n\text{ odd}}}\frac{n-15}{2}q^{n}+
\frac{q^{3}}{1-q^{3}}\sum _{n= 8}^{\infty }q^{2n+1},
\end{eqnarray*}
and by  {\eqref{f3pfg}}, we deduce that $s(n)\geq (n-15)/2$ when $n$ is odd
and $n\geq 17$. Hence, by  {\eqref{f3g}}, we find that
$ f_{3}(n)\geq 0$ for $n\geq 16$ and $f_{3}(n)\geq (n-15)/2$ when
$n$ is odd and $n\geq 17$. Moreover, it can be checked that
$f_{3}(n)\geq 0$ for $2\leq n\leq 15$ and $n\neq 7$. Furthermore,
$f_{3}(3)=f_{3}(5)=0$, $f_{3}(7)=-1$. Thus we complete the proof of  {Theorem~\ref{lem-kernel-2}} when $r=3$.

(4) When $r\geq 4$. Note that
%
\begin{eqnarray}%
\label{equ0sumn0}
\sum _{n=0}^{\infty }f_{r}(n)q^{n}&=&{
\frac{1-q}{(q^{2};q)_{r-1}(q^{2};q)_{r}}}
\\[3pt]
&=&\frac{1}{(q^{2};q)_{r-1}}\sum _{n=0}^{\infty }d_{r+1}(n)q^{n}.
\nonumber
\end{eqnarray}
Since $d_{r+1}(0)=1$ and $d_{r+1}(1)=-1$, we have
%
\begin{eqnarray}
\sum _{n=0}^{\infty }f_{r}(n)q^{n}&=&\frac{1-q}{(q^{2};q)_{r-1}}+
\frac{1}{(q^{2};q)_{r-1}}\sum _{n= 2}^{\infty }d_{r+1}(n)q^{n}
\nonumber
\\[3pt]
&=&\sum _{n=0}^{\infty }d_{r}(n)q^{n}+\frac{1}{(q^{2};q)_{r-1}}\sum _{n=
2}^{\infty }d_{r+1}(n)q^{n}.
\label{pf-main-case-2-3-lemma-1}
\end{eqnarray}
By  {Theorem~\ref{lem-kernel}} (7), we see that when $r\geq 7$,
$d_{r}(n)\geq 0$ for $n\geq 2$ and $d_{r}(2r+7)\geq 1$. Hence, by {\eqref{pf-main-case-2-3-lemma-1}}, we deduce that when $r\geq 7$,
$f_{r}(n)\geq 0$ for $n\geq 2$ and $f_{r}(2r+7)\geq 1$. Applying  {Theorem~\ref{lem-kernel}} (6) and (7), we see that $d_{6}(n)\geq 0$ for
$n\geq 14$ and $d_{7}(n)\geq 0$ for $n\geq 2$, so we derive that
$f_{6}(n)\geq 0$ for $n\geq 14$. It is easy to check that
$f_{6}(n)\geq 0$ for $2\leq n\leq 13$ and $f_{6}(19)\geq 1$. Thus we complete
the proof of  {Theorem~\ref{lem-kernel-2}} when $r\geq 6$.

It remains to show that  {Theorem~\ref{lem-kernel-2}} holds when $r=4$ or
$r=5$. Setting $r=4$ in  {\eqref{pf-main-case-2-3-lemma-1}}, we see that
\begin{eqnarray*}
\sum _{n=0}^{\infty }f_{4}(n)q^{n} = \sum _{n=0}^{\infty }d_{4}(n)q^{n}+
\frac{1}{(q^{2};q)_{3}}\sum _{n= 2}^{\infty }d_{5}(n)q^{n}.
\end{eqnarray*}
From  {Theorem~\ref{lem-kernel}} (5), we see that $d_{5}(n)\geq 0$ for
$n\geq 2$. Moreover, it is easy to see that $d_{5}(2)=1$. Hence
%
\begin{eqnarray}
\sum _{n=0}^{\infty }f_{4}(n)q^{n}&=&\sum _{n= 0}^{\infty }d_{4}(n)q^{n}+
\frac{q^{2}}{(q^{2};q)_{3}}+\frac{1}{(q^{2};q)_{3}}\sum _{n=3}^{\infty }d_{5}(n)q^{n}
\nonumber
\\[3pt]
&=&\sum _{n=0}^{\infty }d_{4}(n)q^{n}+\sum _{n= 2}^{\infty }p_{4}(n-2)q^{n}+
\frac{1}{(q^{2};q)_{3}}\sum _{n=3}^{\infty }d_{5}(n)q^{n}.
\label{pf-lemma-2-3-3a}
\end{eqnarray}
By  {Theorem~\ref{lem-kernel}} (4), we have
\begin{equation*}
d_{4}(n)\geq -\left \lfloor \frac{n+11}{12}\right \rfloor ,
\end{equation*}
and by  {Corollary~\ref{corollary}}, we see that for $n\geq 14$,
\begin{equation*}
p_{4}(n-2)\geq \left \lfloor \frac{n-2}{6}\right \rfloor .
\end{equation*}
Hence for $n\geq 14$,
\begin{equation*}
d_{4}(n)+p_{4}(n-2)\geq \left \lfloor \frac{n-2}{6}\right \rfloor -
\left \lfloor \frac{n+11}{12}\right \rfloor \geq 0.
\end{equation*}
Furthermore, it is routine to check that $d_{4}(n)+p_{4}(n-2)\geq 0$ for
$2\leq n\leq 13$. Note that $d_{5}(n)\geq 0$ for $n\geq 2$, so by  {\eqref{pf-lemma-2-3-3a}}, we conclude that $f_{4}(n)\geq 0$ for
$n\geq 2$. It is easy to check that $f_{4}(15)\geq 1$. Hence  {Theorem~\ref{lem-kernel-2}} is proved when $r=4$.

When $r=5$, by  {\eqref{lem-ker-2-equ}}, we see that
\begin{eqnarray*}
\sum _{n=0}^{\infty }f_{5}(n)q^{n}&=&{\frac{1}{(q^{2};q)_{5}}\cdot
\frac{1-q}{(q^{2};q)_{4}}}
\nonumber
\\[3pt]
&=&\frac{1}{(q^{2};q)_{5}}\sum _{n=0}^{\infty }d_{5}(n)q^{n}.
\end{eqnarray*}
Note that $d_{5}(0)=1$ and $d_{5}(1)=-1$, so we have
%
\begin{eqnarray}%
\label{equ0sumn1}
\sum _{n=0}^{\infty }f_{5}(n)q^{n}&=&\frac{1-q}{(q^{2};q)_{5}}+
\frac{1}{(q^{2};q)_{5}}\sum _{n=2}^{\infty }d_{5}(n)q^{n}
\nonumber
\\[3pt]
&=&\sum _{n=0}^{\infty }d_{6}(n)q^{n}+\frac{1}{(q^{2};q)_{5}}\sum _{n=2}^{\infty }d_{5}(n)q^{n}.
\end{eqnarray}
Since $d_{6}(n)\geq 0$ for $n\geq 14$ and $d_{5}(n)\geq 0$ for
$n\geq 2$, and by  {\eqref{equ0sumn1}}, we derive that $f_{5}(n)\geq 0$ for
$n\geq 14$. It is trivial to check that $f_{5}(n)\geq 0$ for
$2\leq n\leq 13$ and $f_{5}(17)\geq 1$. Hence we arrive at
$f_{5}(n)\geq 0$ for $n\geq 2$ and $f_{5}(17)\geq 1$. Thus, we complete
the proof of  {Theorem~\ref{lem-kernel-2}}.\end{proof}

\section{On $M(m,n)\geq M(m,n-1)$}
\label{sec5}

In this section, we will give a proof of  {Theorem~\ref{main-un-2-n}} by means
of  {Theorem~\ref{lem-kernel}}.

\begin{proof}[Proof of  {Theorem~\ref{main-un-2-n}}] When $m=0$, by  {\eqref{thm-2-equ-n}}, we see that
\begin{eqnarray*}
\sum _{n=0}^{\infty }\left (M(0,n)-M(0,n-1)\right )q^{n} &=&1-2q+q^{2}+q^{3}+
\sum _{k=2}^{\infty
}\frac{q^{k^{2}+2k}}{(q^{2};q)_{k-1}(q^{2};q)_{k-1}}.
\end{eqnarray*}
It yields that $M(0,n)-M(0,n-1)\geq 0$ for $n\geq 2$.

When $m=1$,  {\eqref{thm-2-equ-n}} becomes
\begin{equation*}
\sum _{n=0}^{\infty }\left (M(1,n)-M(1,n-1)\right ) q^{n}=q-q^{2}+
\frac{q^{5}}{1-q^{2}}+\sum _{k=2}^{\infty
}\frac{q^{k^{2}+3k+1}}{(q^{2};q)_{k-1}(q^{2};q)_{k}},
\end{equation*}
which immediately implies that $M(1,n)-M(1,n-1)\geq 0$ for $n\geq 3$.

When $m\geq 2$, by  {\eqref{thm-2-equ-n}}, we have
%
\begin{equation}
\label{thm-2-equ}
M(m,n)-M(m,n-1)\geq d_{m}(n-m)+ p_{m+1}(n-2m-3).
\end{equation}

By  {Theorem~\ref{lem-kernel}} (5)--(7), we see that $d_{m}(n-m)\geq 0$ for
$m \geq 5$ and $n\geq m+2$ except for $(m,n)=(6,13)$ or
$(m,n)=(6,19)$. By the definition of $p_{r}(n)$, we see that
$p_{m+1}(n-2m-3)\geq 0$ for $m\geq 1$. It follows from  {\eqref{thm-2-equ}} that $M(m,n)\geq M(m,n-1)$ for $m \geq 5$ and
$n\geq m+2$ except for $(m,n)=(6,13)$ or $(m,n)=(6,19)$. It is routine
to check that $M(6,19)-M(6,18)\geq 0$ and $M(6,13)-M(6,12)=-1$. Since we
only prove  {\eqref{equ-main-thm-n-un}} holds when $n\ge 14$, we can omit
the case $(m,n)=(6,13)$. Thus  {Theorem~\ref{main-un-2-n}} is verified when
$m\geq 5$.

From  {Theorem~\ref{lem-kernel}} (2)--(4), we see that when
$2\leq m \leq 4$ and $n\geq 4$,
\begin{equation*}
d_{m}(n-m)\geq -\left \lfloor \frac{n+8}{12}\right \rfloor .
\end{equation*}
By  {Corollary~\ref{corollary}}, we derive that when $2\leq m \leq 4$ and
$n\geq 2m+15$,
\begin{equation*}
p_{m+1}(n-2m-3)\geq \left \lfloor \frac{n-2m-3}{6}\right \rfloor .
\end{equation*}
It is easy to check that when $2\leq m \leq 4$ and $n\geq 29$,
\begin{equation*}
\left \lfloor \frac{n-2m-3}{6}\right \rfloor \geq \left \lfloor \frac{n+8}{12}\right \rfloor .
\end{equation*}
So we derive that when $2\leq m \leq 4$ and $n\geq 29$,
\begin{equation*}
M(m,n)-M(m,n-1)\geq d_{m}(n-m)+p_{m+1}(n-2m-3)\geq 0.
\end{equation*}
Moreover, it can be checked that $M(m,n)\geq M(m,n-1)$ when
$2\leq m \leq 4$ and $14\leq n\leq 28$. So  {Theorem~\ref{main-un-2-n}} is
verified when $2\leq m \leq 4$. Thus, we complete the proof of  {Theorem~\ref{main-un-2-n}}. \end{proof}

\section{The generating function of $M(m-1,n)-M(m,n)$}
\label{sec6}

In this section, we will establish three expressions for the generating
function of $M(m-1,n)-M(m,n)$, which play a crucial role in the proof of
 {Theorem~\ref{main-thm-n}}. To this end, we first split the first summation
in  {\eqref{lem-gf-r-g}} into five summations as follows.

\begin{lem}%
\label{lem-m-ge-1-gen}
For $m\geq 1$,
%
\begin{eqnarray}%
\label{lem-m-ge-1-gen-eq}
{\sum _{k=1}^{\infty
}\frac{q^{k(k+m-1)+2k+m-1}}{(q;q)_{k}(q^{2};q)_{k+m-2}}}&=&\frac{q^{2m+2}}{(q^{2};q)_{m-1}}+
\frac{q^{3m+7}}{(1-q^{2})(q^{2};q)_{m-1}}
\nonumber
\\[3pt]
&&+\sum _{k=1}^{\infty
}\frac{q^{k(k+m)+k+m-1}}{(q;q)_{k-1}(q^{2};q)_{k+m-2}}+\sum _{k=3}^{\infty }\frac{q^{k(k+m)+2k+m-1}}{(q^{2};q)_{k-1}(q^{2};q)_{k+m-4}}
\nonumber
\\[3pt]
&&+\sum _{k= 1}^{\infty
}\frac{q^{k(k+m)+2k+m}}{(q;q)_{k}(q^{2};q)_{k+m-2}}+\sum _{k=3}^{\infty }\frac{q^{k(k+m)+3k+2m-3}}{(q^{2};q)_{k-1}(q^{2};q)_{k+m-3}}
\nonumber
\\[3pt]
&&+\sum _{k=2}^{\infty
}\frac{q^{k(k+m)+3k+2m-2}}{(q^{2};q)_{k-1}(q^{2};q)_{k+m-2}}.
\label{main-p-final}
\end{eqnarray}
\end{lem}

\begin{proof} It is clear that when $m\geq 1$,
%
\begin{equation}
\label{equ-new-sum-n-0-mm-1}
{\sum _{k=1}^{\infty
}\frac{q^{k(k+m-1)+2k+m-1}}{(q;q)_{k}(q^{2};q)_{k+m-2}}=}\sum _{k=1}^{\infty }\frac{q^{k(k+m)+k+m-1}}{(q;q)_{k-1}(q^{2};q)_{k+m-2}}\cdot
\frac{1}{1-q^{k}}
\end{equation}
Obviously, when $k\geq 1$,
%
\begin{equation}
\label{equ-new-sum-n-0-mm-2}
\frac{1}{1-q^{k}}=1+\frac{q^{k}(1-q)}{1-q^{k}}+
\frac{q^{k+1}}{1-q^{k}}.
\end{equation}
Substituting  {\eqref{equ-new-sum-n-0-mm-2}} into  {\eqref{equ-new-sum-n-0-mm-1}}, we deduce that
%
\begin{eqnarray}%
\label{equ-new-sum-n-0-mm-3}
{\sum _{k=1}^{\infty
}\frac{q^{k(k+m-1)+2k+m-1}}{(q;q)_{k}(q^{2};q)_{k+m-2}}} &=&\sum _{k=1}^{\infty }\frac{q^{k(k+m)+k+m-1}}{(q;q)_{k-1}(q^{2};q)_{k+m-2}}+\sum _{k=1}^{\infty }\frac{q^{k(k+m)+2k+m-1}}{(q^{2};q)_{k-1}(q^{2};q)_{k+m-2}}
\nonumber
\\
&&+\sum _{k=1}^{\infty
}\frac{q^{k(k+m)+2k+m}}{(q;q)_{k}(q^{2};q)_{k+m-2}}.
\end{eqnarray}

Notice that
%
\begin{eqnarray}%
\label{equ-new-sum-n-0-mm-4}
&&\sum _{k=1}^{\infty
}\frac{q^{k(k+m)+2k+m-1}}{(q^{2};q)_{k-1}(q^{2};q)_{k+m-2}}
\nonumber
\\[3pt]
&=&\frac{q^{2m+2}}{(q^{2};q)_{m-1}}+ \sum _{k=2}^{\infty
}\frac{q^{k(k+m)+2k+m-1}}{(q^{2};q)_{k-1}(q^{2};q)_{k+m-2}}
\nonumber
\\[3pt]
&=&\frac{q^{2m+2}}{(q^{2};q)_{m-1}}+\sum _{k=2}^{\infty
}\frac{q^{k(k+m)+2k+m-1}}{(q^{2};q)_{k-1}(q^{2};q)_{k+m-3}} \left (1+
\frac{q^{k+m-1}}{1-q^{k+m-1}}\right )
\nonumber
\\[3pt]
&=&\frac{q^{2m+2}}{(q^{2};q)_{m-1}}+\sum _{k=2}^{\infty
}\frac{q^{k(k+m)+3k+2m-2}}{(q^{2};q)_{k-1}(q^{2};q)_{k+m-2}}+\sum _{k= 2}^{\infty }\frac{q^{k(k+m)+2k+m-1}}{(q^{2};q)_{k-1}(q^{2};q)_{k+m-3}}
\nonumber
\\[3pt]
&=&\frac{q^{2m+2}}{(q^{2};q)_{m-1}}+\sum _{k=2}^{\infty
}\frac{q^{k(k+m)+3k+2m-2}}{(q^{2};q)_{k-1}(q^{2};q)_{k+m-2}} +
\frac{q^{3m+7}}{(1-q^{2})(q^{2};q)_{m-1}}
\nonumber
\\[3pt]
&&+\sum _{k=3}^{\infty
}\frac{q^{k(k+m)+2k+m-1}}{(q^{2};q)_{k-1}(q^{2};q)_{k+m-4}} \left (1+
\frac{q^{k+m-2}}{1-q^{k+m-2}}\right )
\nonumber
\\[3pt]
&=&\frac{q^{2m+2}}{(q^{2};q)_{m-1}}+
\frac{q^{3m+7}}{(1-q^{2})(q^{2};q)_{m-1}}+\sum _{k=2}^{\infty
}\frac{q^{k(k+m)+3k+2m-2}}{(q^{2};q)_{k-1}(q^{2};q)_{k+m-2}}
\nonumber\\
&&+\sum _{k= 3}^{\infty }\frac{q^{k(k+m)+3k+2m-3}}{(q^{2};q)_{k-1}(q^{2};q)_{k+m-3}}
\nonumber
\\[3pt]
&&+\sum _{k=3}^{\infty
}\frac{q^{k(k+m)+2k+m-1}}{(q^{2};q)_{k-1}(q^{2};q)_{k+m-4}}.
\label{main-p-3-4-5}
\end{eqnarray}
Substituting  {\eqref{main-p-3-4-5}} into  {\eqref{equ-new-sum-n-0-mm-3}}, we
are led to  {\eqref{main-p-final}}, and hence  {Lemma~\ref{lem-m-ge-1-gen}} follows.\end{proof}

The second summation in  {\eqref{lem-gf-r-g}} can be split into the following
five summations.

\begin{lem}%
\label{lem-m-ge-2-gen}%
For $m\geq 1$,
%
\begin{eqnarray}
{\sum _{k=1}^{\infty
}\frac{q^{k(k+m)+2k+m}}{(q;q)_{k}(q^{2};q)_{k+m-1}}}&=&\sum _{k=1}^{\infty }\frac{q^{k(k+m)+2k+m}}{(q;q)_{k}(q^{2};q)_{k+m-2}}+\sum _{k=1}^{\infty }\frac{q^{k(k+m)+3k+2m}}{(q^{2};q)_{k-1}(q^{2};q)_{k+m-2}}
\nonumber
\\[3pt]
&&+\sum _{k=1}^{\infty
}\frac{q^{k(k+m)+3k+2m+1}}{(q;q)_{k}(q^{2};q)_{k+m-1}}+\sum _{k=1}^{\infty }\frac{q^{k(k+m)+4k+3m}}{(q^{2};q)_{k-1}(q^{2};q)_{k+m-2}}
\nonumber
\\[3pt]
&&+\sum _{k=1}^{\infty
}\frac{q^{k(k+m)+5k+4m}}{(q^{2};q)_{k-1}(q^{2};q)_{k+m-1}}.
\label{main-s-final}
\end{eqnarray}
\end{lem}

\begin{proof} Clearly,
%
\begin{equation}
\label{main-s1}
{\sum _{k=1}^{\infty
}\frac{q^{k(k+m)+2k+m}}{(q;q)_{k}(q^{2};q)_{k+m-1}}=} \sum _{k=1}^{\infty }\frac{q^{k(k+m)+2k+m}}{(q;q)_{k}(q^{2};q)_{k+m-2}}\cdot
\frac{1}{1-q^{k+m}}.
\end{equation}
Moreover, one can easily check that the following identity holds:
%
\begin{equation}
\label{equ-1-frac-q-km}
\frac{1}{1-q^{k+m}}=1+q^{k+m}(1-q)+\frac{q^{k+m+1}}{1-q^{k+m}}+ q^{2k+2m}(1-q)+
\frac{q^{3k+3m}(1-q)}{1-q^{k+m}}.
\end{equation}

Substituting  {\eqref{equ-1-frac-q-km}} into  {\eqref{main-s1}}, we obtain  {\eqref{main-s-final}}. This completes the proof.\end{proof}

By  {Lemma~\ref{lem-m-ge-1-gen}} and  {Lemma~\ref{lem-m-ge-2-gen}}, we obtain
the first expression of the generating function of $M(m-1,n)-M(m,n)$ when
$m\geq 1$.

\begin{thm}%
\label{lem-4-1}
For $m\geq 1$,
%
\begin{eqnarray}
&&\sum _{n=0}^{\infty }\left (M(m-1,n)-M(m,n)\right )q^{n}
\nonumber
\\[3pt]
&=&\frac{q^{m-1}(1-q)}{(q;q)_{m-1}}-\frac{q^{m}}{(q^{2};q)_{m-1}} +
\frac{q^{2m+1}}{(q^{2};q)_{m-1}}+\frac{q^{2m+2}}{(q^{2};q)_{m-1}}
\nonumber
\\[3pt]
&&-\frac{q^{3m+4}}{(q^{2};q)_{m-1}}+
\frac{q^{3m+7}}{(1-q^{2})(q^{2};q)_{m-1}}-
\frac{q^{5m+6}}{(q^{2};q)_{m}}
\nonumber
\\[3pt]
&&+\sum _{k=3}^{\infty
}\frac{q^{k(k+m)+2k+m-1}}{(q^{2};q)_{k-1}(q^{2};q)_{k+m-4}}+\sum _{k=2}^{\infty }\frac{q^{k(k+m)+3k+2m-2}}{(q^{3};q)_{k-2}(q^{2};q)_{k+m-2}}
\nonumber
\\[3pt]
&&-\sum _{k=1}^{\infty
}\frac{q^{k(k+m)+4k+3m}}{(q^{2};q)_{k-1}(q^{2};q)_{k+m-2}}+\sum _{k=2}^{\infty
}\frac{q^{k(k+m)+5k+3m+1}(1-q^{m-1})}{(q^{2};q)_{k}(q^{2};q)_{k+m-1}}.
\label{lem-4-1-main-final}
\end{eqnarray}
\end{thm}

\begin{proof} Substituting  {\eqref{main-p-final}} and  {\eqref{main-s-final}} into {\eqref{lem-gf-r-g}}, and by simplification, we
get
%
\begin{eqnarray}%
\label{main-uncr-1}
&&\sum _{n=0}^{\infty }\left (M(m-1,n)-M(m,n)\right )q^{n}
\\[3pt]
&=&\left (\sum _{k=2}^{\infty
}\frac{q^{k(k+m)+3k+2m-2}}{(q^{2};q)_{k-1}(q^{2};q)_{k+m-2}} -\sum _{k=1}^{\infty }\frac{q^{k(k+m)+3k+2m}}{(q^{2};q)_{k-1}(q^{2};q)_{k+m-2}}
\right )
\nonumber
\label{m-4-1-case-3}
\\[3pt]
&&+\left (\sum _{k=3}^{\infty
}\frac{q^{k(k+m)+3k+2m-3}}{(q^{2};q)_{k-1}(q^{2};q)_{k+m-3}}-\sum _{k=1}^{\infty }\frac{q^{k(k+m)+5k+4m}}{(q^{2};q)_{k-1}(q^{2};q)_{k+m-1}}
\right )
\nonumber
\label{m-4-1-case-4}
\\[3pt]
&&+ \left (\sum _{k=1}^{\infty
}\frac{q^{k(k+m)+k+m-1}}{(q;q)_{k-1}(q^{2};q)_{k+m-2}}-\sum _{k=1}^{\infty }\frac{q^{k(k+m)+3k+2m+1}}{(q;q)_{k}(q^{2};q)_{k+m-1}}\right )
\nonumber
\\[3pt]
&&{+\sum _{k=3}^{\infty
}\frac{q^{k(k+m)+2k+m-1}}{(q^{2};q)_{k-1}(q^{2};q)_{k+m-4}}-\sum _{k=1}^{\infty }\frac{q^{k(k+m)+4k+3m}}{(q^{2};q)_{k-1}(q^{2};q)_{k+m-2}}
\label{main-dis-1}
}
\nonumber
\\[3pt]
&&+ \frac{q^{m-1}(1-q)}{(q;q)_{m-1}}-\frac{q^{m}}{(q^{2};q)_{m-1}} +
\frac{q^{2m+2}}{(q^{2};q)_{m-1}}+
\frac{q^{3m+7}}{(1-q^{2})(q^{2};q)_{m-1}}.
\nonumber
\label{m-4-1-case-5}
\end{eqnarray}
Observe that
%
\begin{eqnarray}%
\label{case-3-tt-1}
&&\sum _{k=2}^{\infty
}\frac{q^{k(k+m)+3k+2m-2}}{(q^{2};q)_{k-1}(q^{2};q)_{k+m-2}}-\sum _{k=1}^{\infty }\frac{q^{k(k+m)+3k+2m}}{(q^{2};q)_{k-1}(q^{2};q)_{k+m-2}}
\nonumber
\\[3pt]
&=&-\frac{q^{3m+4}}{(q^{2};q)_{m-1}}+\sum _{k=2}^{\infty
}\frac{q^{k(k+m)+3k+2m-2}(1-q^{2})}{(q^{2};q)_{k-1}(q^{2};q)_{k+m-2}}
\nonumber
\\[3pt]
&=&-\frac{q^{3m+4}}{(q^{2};q)_{m-1}}+\sum _{k=2}^{\infty
}\frac{q^{k(k+m)+3k+2m-2} }{(q^{3};q)_{k-2}(q^{2};q)_{k+m-2}},
\end{eqnarray}
and
%
\begin{eqnarray}
&&\sum _{k=3}^{\infty
}\frac{q^{k(k+m)+3k+2m-3}}{(q^{2};q)_{k-1}(q^{2};q)_{k+m-3}}-\sum _{k=1}^{\infty }\frac{q^{k(k+m)+5k+4m}}{(q^{2};q)_{k-1}(q^{2};q)_{k+m-1}}
\nonumber
\\[3pt]
&=& \sum _{k=2}^{\infty
}\frac{q^{k(k+m)+5k+3m+1}}{(q^{2};q)_{k}(q^{2};q)_{k+m-2}}- \sum _{k=1}^{\infty }\frac{q^{k(k+m)+5k+4m}}{(q^{2};q)_{k-1}(q^{2};q)_{k+m-1}}
\nonumber
\\[3pt]
&=&-\frac{q^{5m+6}}{(q^{2};q)_{m}}+\sum _{k=2}^{\infty
}\frac{q^{k(k+m)+5k+3m+1}}{(q^{2};q)_{k}(q^{2};q)_{k+m-1}}\left ((1-q^{k+m})-q^{m-1}(1-q^{k+1})
\right )
\nonumber
\\[3pt]
&=&-\frac{q^{5m+6}}{(q^{2};q)_{m}}+\sum _{k=2}^{\infty
}\frac{q^{k(k+m)+5k+3m+1}(1-q^{m-1})}{(q^{2};q)_{k}(q^{2};q)_{k+m-1}}.
\label{case-1-4}
\end{eqnarray}
Moreover, it is easy to see that
%
\begin{equation}
\label{lem-4-1-main-final=abc}
\sum _{k=1}^{\infty
}\frac{q^{k(k+m)+k+m-1}}{(q;q)_{k-1}(q^{2};q)_{k+m-2}}-\sum _{k=1}^{\infty }\frac{q^{k(k+m)+3k+2m+1}}{(q;q)_{k}(q^{2};q)_{k+m-1}}=
\frac{q^{2m+1}}{(q^{2};q)_{m-1}}.
\end{equation}
We then obtain  {\eqref{lem-4-1-main-final}} upon substituting  {\eqref{case-3-tt-1}}, {\eqref{case-1-4}} and  {\eqref{lem-4-1-main-final=abc}} into  {\eqref{main-uncr-1}}. This completes
the proof.\end{proof}

When $m\ge 2$, we find that the generating function of
$M(m-1,n)-M(m,n)$ in  {Theorem~\ref{lem-4-1}} can be further simplified as stated below.

\begin{thm}%
\label{lem-4-2}
For $m\geq 2$,
%
\begin{eqnarray}
&&\sum _{n=0}^{\infty }\left (M(m-1,n)-M(m,n)\right )q^{n}
\nonumber
\\[3pt]
&=&\frac{q^{m-1}}{(q^{2};q)_{m-2}} -\frac{q^{m}}{(q^{2};q)_{m-2}}-
\frac{q^{2m} }{(q^{3};q)_{m-2}}+\frac{q^{2m+1}}{(q^{2};q)_{m-1}}-
\frac{q^{3m+4}}{(q^{2};q)_{m-1}}
\nonumber
\\[3pt]
&&+\sum _{k=2}^{\infty
}\frac{q^{k(k+m)+3k+2m-2}}{(q^{3};q)_{k-2}(q^{2};q)_{k+m-2}}+\sum _{k=1}^{\infty
}\frac{q^{k(k+m)+4k+2m+2}(1-q^{m-2})}{(q^{2};q)_{k}(q^{2};q)_{k+m-2}}
\nonumber
\\[3pt]
&&+\sum _{k=1}^{\infty
}\frac{q^{k(k+m)+5k+3m+1}(1-q^{m-1})}{(q^{2};q)_{k}(q^{2};q)_{k+m-1}}.
\label{lem-4-2-main-final}
\end{eqnarray}
\end{thm}

\begin{proof} It is trivial to verify that when $m\ge 2$,
%
\begin{equation}
\frac{q^{m}}{(q^{2};q)_{m-1}}=\frac{q^{m}}{(q^{2};q)_{m-2}}+
\frac{q^{2m} }{(q^{3};q)_{m-2}}+\frac{q^{2m+2} }{(q^{2};q)_{m-1}}.
\end{equation}
Hence, by  {Theorem~\ref{lem-4-1}}, it suffices to show that
%
\begin{eqnarray}
&&\frac{q^{3m+7}}{(1-q^{2})(q^{2};q)_{m-1}}-
\frac{q^{5m+6}}{(q^{2};q)_{m}}+\sum _{k=3}^{\infty
}\frac{q^{k(k+m)+2k+m-1}}{(q^{2};q)_{k-1}(q^{2};q)_{k+m-4}}
\nonumber
\\[3pt]
&&-\sum _{k=1}^{\infty
}\frac{q^{k(k+m)+4k+3m}}{(q^{2};q)_{k-1}(q^{2};q)_{k+m-2}}+\sum _{k=2}^{\infty
}\frac{q^{k(k+m)+5k+3m+1}(1-q^{m-1})}{(q^{2};q)_{k}(q^{2};q)_{k+m-1}}
\nonumber
\\[7pt]
&=&\sum _{k=1}^{\infty
}\frac{q^{k(k+m)+4k+2m+2}(1-q^{m-2})}{(q^{2};q)_{k}(q^{2};q)_{k+m-2}}+
\sum _{k=1}^{\infty
}\frac{q^{k(k+m)+5k+3m+1}(1-q^{m-1})}{(q^{2};q)_{k}(q^{2};q)_{k+m-1}}.
\label{lem-4-2-case-7}
\end{eqnarray}

First, observe that
%
\begin{eqnarray}%
\label{lem-4-2-case-7a}
&&\sum _{k=3}^{\infty
}\frac{q^{k(k+m)+2k+m-1}}{(q^{2};q)_{k-1}(q^{2};q)_{k+m-4}}-\sum _{k=1}^{\infty }\frac{q^{k(k+m)+4k+3m}}{(q^{2};q)_{k-1}(q^{2};q)_{k+m-2}}
\nonumber
\\[3pt]
&=&\sum _{k=2}^{\infty
}\frac{q^{k(k+m)+4k+2m+2}}{(q^{2};q)_{k}(q^{2};q)_{k+m-3}}-\sum _{k=1}^{\infty }\frac{q^{k(k+m)+4k+3m}}{(q^{2};q)_{k-1}(q^{2};q)_{k+m-2}}
\nonumber
\\[3pt]
&=&-\frac{q^{4m+5}}{(q^{2};q)_{m-1}}+\sum _{k=2}^{\infty
}\frac{q^{k(k+m)+4k+2m+2}(1-q^{m-2})}{(q^{2};q)_{k}(q^{2};q)_{k+m-2}}.
\end{eqnarray}
On the other hand, we find that when $m\geq 2$,
%
\begin{eqnarray}%
\label{lem-4-2-case-7b}
&&\quad \frac{q^{3m+7}}{(1-q^{2})(q^{2};q)_{m-1}}-
\frac{q^{4m+5}}{(q^{2};q)_{m-1}}-\frac{q^{5m+6}}{(q^{2};q)_{m}}
\nonumber
\\[3pt]
&&=\frac{q^{3m+7}}{(1-q^{2})(q^{2};q)_{m-2}}\left (1+
\frac{q^{m}}{1-q^{m}}\right )- \frac{q^{4m+5}}{(q^{2};q)_{m-1}}-
\frac{q^{5m+6}}{(q^{2};q)_{m}}
\nonumber
\\[3pt]
&&=\frac{q^{3m+7}}{(1-q^{2})(q^{2};q)_{m-2}}-
\frac{q^{4m+5}}{(q^{2};q)_{m-1}}+
\frac{q^{4m+7}}{(1-q^{2})(q^{2};q)_{m-1}}-
\frac{q^{5m+6}}{(q^{2};q)_{m}}
\nonumber
\\[3pt]
&&= \frac{q^{3m+7}(1-q^{m-2})}{(1-q^{2})(q^{2};q)_{m-1}}+
\frac{q^{4m+7}(1-q^{m-1})}{(1-q^{2})(q^{2};q)_{m}}.
\end{eqnarray}
Substituting  {\eqref{lem-4-2-case-7a}} and  {\eqref{lem-4-2-case-7b}} into the
left-hand side of  {\eqref{lem-4-2-case-7}}, we obtain the right-hand side
of  {\eqref{lem-4-2-case-7}}. This completes the proof of  {Theorem~\ref{lem-4-1}}.\end{proof}

When $m\geq 3$, we could further simplify the generating function of
$M(m-1,n)-M(m,n)$ in  {Theorem~\ref{lem-4-2}} to obtain the following expression.

\begin{thm}%
\label{lem-4}
For $m\geq 3$,
%
\begin{eqnarray}
&&\sum _{n=0}^{\infty }\left (M(m-1,n)-M(m,n)\right )q^{n}
\nonumber
\\[3pt]
&=&-q^{2m}+q^{2m+1}+q^{3m+1}+\frac{q^{m-1}}{(q^{2};q)_{m-2}} -
\frac{q^{m}}{(q^{2};q)_{m-2}}+
\frac{q^{2m+5}}{(q^{2};q)_{m-3}(1-q^{m})}
\nonumber
\\[3pt]
&&+ \sum _{k=3}^{m}\frac{q^{2k+2m+1}}{(q^{k};q)_{m-k+1}}+\sum _{k=2}^{\infty }\frac{q^{k(k+m)+3k+2m-2}}{(q^{3};q)_{k-2}(q^{2};q)_{k+m-2}}+
\sum _{k=1}^{\infty
}\frac{q^{k(k+m)+4k+2m+2}(1-q^{m-2})}{(q^{2};q)_{k}(q^{2};q)_{k+m-2}}
\nonumber
\\[3pt]
&&+\sum _{k=1}^{\infty
}\frac{q^{k(k+m)+5k+3m+1}}{(q^{2};q)_{k}(q^{2};q)_{m-3}(q^{m};q)_{k+1}}.
\label{main-final}
\end{eqnarray}
\end{thm}

\begin{proof} From  {Theorem~\ref{lem-4-2}}, it suffices to show that when
$m\geq 3$,
%
\begin{eqnarray}
&&\frac{q^{2m+1}}{(q^{2};q)_{m-1}}-\frac{q^{2m} }{(q^{3};q)_{m-2}}-
\frac{q^{3m+4}}{(q^{2};q)_{m-1}}
\nonumber
\\[3pt]
&=& -q^{2m}+q^{2m+1}+q^{3m+1}+
\frac{q^{2m+5}}{(q^{2};q)_{m-3}(1-q^{m})} + \sum _{k=3}^{m}
\frac{q^{2k+2m+1}}{(q^{k};q)_{m-k+1}}
\label{equ-2-lem-4-2-final-1}
.
\end{eqnarray}
In light of  {\eqref{lem-snk-gf-3}}, we see that
%
\begin{eqnarray}%
\label{equ-3-gf-first}
\frac{q^{2m+1}}{(q^{2};q)_{m-1}}&=& q^{2m+1}+\sum _{k=2}^{m}
\frac{q^{2m+1+k}}{(q^{k};q)_{m-k+1}}
\nonumber
\\[3pt]
&=&q^{2m+1}+\sum _{k=2}^{m}
\frac{q^{2m+1+k}(1-q^{k}+q^{k})}{(q^{k};q)_{m-k+1}}
\nonumber
\\[3pt]
&=&q^{2m+1}+\sum _{k=2}^{m}\frac{q^{k+2m+1}}{(q^{k+1};q)_{m-k}}+
\sum _{k=2}^{m}\frac{q^{2k+2m+1}}{(q^{k};q)_{m-k+1}}.
\label{temp-1}
\end{eqnarray}
Using the same argument as in the proof of  {\eqref{lem-snk-gf-3}}, we deduce
that for $m\geq 3$,
%
\begin{equation}
\label{equ-3-gf-sec}
\frac{q^{2m}}{(q^{3};q)_{m-2}}=q^{2m}+\sum _{k=3}^{m}
\frac{q^{k+2m}}{(q^{k};q)_{m-k+1}} =q^{2m}+\sum _{k=2}^{m-1}
\frac{q^{k+2m+1}}{(q^{k+1};q)_{m-k}}.
\end{equation}
Substituting  {\eqref{equ-3-gf-first}} and  {\eqref{equ-3-gf-sec}} into the left-hand
side of  {\eqref{equ-2-lem-4-2-final-1}}, we obtain
\begin{eqnarray*}
&&\frac{q^{2m+1}}{(q^{2};q)_{m-1}}-\frac{q^{2m}}{(q^{3};q)_{m-2}}-
\frac{q^{3m+4}}{(q^{2};q)_{m-1}}
\nonumber
\\[3pt]
&=&-q^{2m}+q^{2m+1}+q^{3m+1}+ \sum _{k=3}^{m}
\frac{q^{2k+2m+1}}{(q^{k};q)_{m-k+1}}+
\frac{q^{2m+5}}{(q^{2};q)_{m-1}}-\frac{q^{3m+4}}{(q^{2};q)_{m-1}}
\nonumber
\\[3pt]
&=&-q^{2m}+q^{2m+1}+q^{3m+1}+ \sum _{k=3}^{m}
\frac{q^{2k+2m+1}}{(q^{k};q)_{m-k+1}}+
\frac{q^{2m+5}}{(q^{2};q)_{m-3}(1-q^{m})},
\end{eqnarray*}
which is equal to the right-hand side of  {\eqref{equ-2-lem-4-2-final-1}}. Thus, we complete the proof of {Theorem~\ref{lem-4}}.\end{proof}

\section{On $M(1,n)\geq M(2,n)$}
\label{sec7}

In the following three sections, we will give a proof of  {Theorem~\ref{main-thm-n}}. In this section, we will show that  {Theorem~\ref{main-thm-n}} holds when $m=2$. In Section~\ref{sec8}, we will prove that {Theorem~\ref{main-thm-n}} holds when $m\geq 3$. Section~\ref{sec9} is devoted to the proof of {Theorem~\ref{main-thm-n}}  when $m=1$. As  stated in Section~\ref{sec2}, the
proof of  {Theorem~\ref{main-thm-n}} when $m=1$ is the most complicated, so
we put the proof of the case $m=1$ at the end of the proof of the whole
theorem.

\begin{proof}[Proof of  {Theorem~\ref{main-thm-n}} for $m=2$] Setting $m=2$ in {Theorem~\ref{lem-4-2}}, we have
%
\begin{eqnarray}
&&\sum _{n=0}^{\infty }\left (M(1,n)-M(2,n)\right )q^{n}
\nonumber
\\[3pt]
&=&q-q^{2}-q^{4}+\frac{q^{5}}{1-q^{2}}-\frac{q^{10}}{1-q^{2}} +\sum _{k=2}^{\infty }\frac{q^{k^{2}+5k+2}}{(q^{3};q)_{k-2}(q^{2};q)_{k}}
\nonumber
\\[3pt]
&&+\sum _{k=1}^{\infty
}\frac{q^{k^{2}+7k+7}(1-q)}{(q^{2};q)_{k}(q^{2};q)_{k+1}}.
\label{main-m-2-dis-76}
\end{eqnarray}
Observe that
%
\begin{eqnarray}
&&\sum _{k=2}^{\infty
}\frac{q^{k^{2}+5k+2}}{(q^{3};q)_{k-2}(q^{2};q)_{k}}-
\frac{q^{10}}{1-q^{2}}
\nonumber
\\[3pt]
&=&\sum _{k=3}^{\infty
}\frac{q^{k^{2}+5k+2}}{(q^{3};q)_{k-2}(q^{2};q)_{k}}+
\frac{q^{16}}{(1-q^{2})(1-q^{3})}-\frac{q^{10}}{1-q^{2}}
\nonumber
\\[3pt]
&=&\sum _{k=3}^{\infty
}\frac{q^{k^{2}+5k+2}}{(q^{3};q)_{k-2}(q^{2};q)_{k}}+
\frac{q^{19}}{(1-q^{2})(1-q^{3})} -q^{10}-q^{12}-q^{14}.
\label{main-m-2-1}
\end{eqnarray}
Define
%
\begin{equation}
\label{main-m-3-claim}
\sum _{n=0}^{\infty }T_{2}(n)q^{n}:=\sum _{k=1}^{\infty
}\frac{q^{k^{2}+7k+7}(1-q)}{(q^{2};q)_{k}(q^{2};q)_{k+1}},
\end{equation}
and by  {\eqref{main-m-2-dis-76}} and  {\eqref{main-m-2-1}}, we find that for
$n\geq 15$,
%
\begin{equation}
M(1,n)-M(2,n)\geq T_{2}(n).
\end{equation}
Hence it suffices to show that $T_{2}(n)\geq 0$ when $n\geq 44$.

By  {\eqref{lem-ker-2-equ}} and  {\eqref{main-m-3-claim}}, we find that
$T_{2}(n)$ can be expressed in terms of $f_{r}(n)$ as follows.
%
\begin{eqnarray}%
\label{main-m-3-claim-pf}
\sum _{n=0}^{\infty }T_{2}(n)q^{n}=\sum _{k=1}^{\infty
}\frac{q^{k^{2}+7k+7}(1-q)}{(q^{2};q)_{k}(q^{2};q)_{k+1}}= \sum _{k=1}^{\infty }q^{k^{2}+7k+7}\sum _{n=0}^{\infty }f_{k+1}(n)q^{n}.
\end{eqnarray}
Define
\begin{eqnarray*}
\sum _{n=0}^{\infty }R(n)q^{n}&:=&q^{15}\sum _{n=0}^{\infty }f_{2}(n)q^{n}+q^{25}
\sum _{n=0}^{\infty }f_{3}(n)q^{n},
\\[3pt]
\sum _{n=0}^{\infty }S(n)q^{n}&:=&\sum _{k=3}^{\infty }q^{k^{2}+7k+7}
\sum _{n=0}^{\infty }f_{k+1}(n)q^{n}.
\end{eqnarray*}
By  {\eqref{main-m-3-claim-pf}}, we find that for $n\geq 0$,
%
\begin{equation}
\label{main-m-3-claim-pf-re}
T_{2}(n)=R(n)+S(n).
\end{equation}
We will investigate the nonnegativity of $R(n)$ and $S(n)$ respectively.

By  {Theorem~\ref{lem-kernel-2}} (2), we see that
%
\begin{eqnarray}%
\label{main-m-3-1}
q^{15}\sum _{n=0}^{\infty }f_{2}(n)q^{n} &=&q^{15}\left (\sum _{m= 0}^{\infty }f_{2}(2m)q^{2m}-\sum _{m=0}^{\infty }\left \lceil \frac{2m+1}{6}
\right \rceil q^{2m+1}\right )
\nonumber
\\[3pt]
&=&\sum _{m=7}^{\infty }f_{2}(2m-14)q^{2m+1}-\sum _{m= 8}^{\infty }\left
\lceil \frac{2m-15}{6}\right \rceil q^{2m}.
\end{eqnarray}
From  {Theorem~\ref{lem-kernel-2}} (3), we have
%
\begin{align}
\label{main-m-3-2}
q^{25}\sum _{n=0}^{\infty }f_{3}(n)q^{n} =&q^{25}\left (1-q-q^{7}+\sum _{m=4}^{\infty }f_{3}(2m+1)q^{2m+1}+\sum _{m= 1}^{\infty }f_{3}(2m)q^{2m}\right )
\nonumber
\\[3pt]
=&q^{25}-q^{26}-q^{32}+\sum _{m=17}^{\infty }f_{3}(2m-25)q^{2m}+\sum _{m=13}^{\infty }f_{3}(2m-24)q^{2m+1}.
\end{align}
Combining  {\eqref{main-m-3-1}} and  {\eqref{main-m-3-2}}, we find that
\begin{eqnarray*}
\sum _{n=15}^{\infty }R(n)q^{n}&=&q^{25}-q^{26}-q^{32} -\sum _{m=8}^{20}
\left \lceil \frac{2m-15}{6}\right \rceil q^{2m}
\nonumber
\\[3pt]
&&+\sum _{m=7}^{\infty }f_{2}(2m-14)q^{2m+1}+\sum _{m=17}^{20}f_{3}(2m-25)q^{2m}
\nonumber\\
&&+
\sum _{m= 13}^{\infty }f_{3}(2m-24)q^{2m+1}
\nonumber
\\[3pt]
&&+\sum _{m= 21}^{\infty }\left (f_{3}(2m-25)- \left \lceil \frac{2m-15}{6}\right \rceil \right ) q^{2m} .
\end{eqnarray*}
By  {Theorem~\ref{lem-kernel-2}} (2) and (3), we see that $f_{2}(2m-14)
\geq 0$ when $m\geq 7$ and $f_{3}(2m-24)\ge 0$ when $m\ge 13$. It yields
that when $m\geq 7$,
%
\begin{equation}
\label{main-m-3-12a}
R(2m+1)\geq 0.
\end{equation}
Using  {Theorem~\ref{lem-kernel-2}} (3), we see that
$f_{3}(2m-25)\geq m-20$ for $m\geq 21$. It follows that when
$m\geq 27$,
\begin{equation*}
f_{3}(2m-25)- \left \lceil \frac{2m-15}{6}\right \rceil \geq m-20-
\left \lceil \frac{2m-15}{6}\right \rceil \geq 0.
\end{equation*}
Hence we derive that when $m\geq 27$,
%
\begin{equation}
\label{main-m-3-12b}
R(2m)\geq 0.
\end{equation}
Combining  {\eqref{main-m-3-12a}} and  {\eqref{main-m-3-12b}}, we derive that
when $n\geq 54$,
\begin{equation*}
R(n)\geq 0.
\end{equation*}
It can be checked that $R(n)\geq 0$ for $44\leq n\leq 53$. Thus we show
that $R(n)\geq 0$ for $n\geq 44$.

We proceed to investigate the nonnegativity of $S(n)$. Observe that
%
\begin{eqnarray}
\sum _{n=0}^{\infty }S(n)q^{n}&=& {\sum _{k=3}^{\infty }q^{k^{2}+7k+7}
\sum _{n=0}^{\infty }f_{k+1}(n)q^{n}}
\nonumber
\\[3pt]
&=&\sum _{k=3}^{\infty }q^{k^{2}+7k+7}\left (1-q +f_{k+1}(2k+9)q^{2k+9}+
\sum _{{n\geq 2 }\atop{n\neq 2k+9}}f_{k+1}(n)q^{n}\right )
\nonumber
\\[3pt]
&=&\sum _{k=3}^{\infty }q^{k^{2}+7k+7}\left (1+\sum _{{n\geq 2 }\atop{n
\neq 2k+9}}f_{k+1}(n)q^{n}\right )
\nonumber
\\[3pt]
&&+\sum _{k=3}^{\infty }f_{k+1}(2k+9)q^{k^{2}+9k+16}-\sum _{k=3}^{\infty }q^{k^{2}+7k+8}.
\label{main-m-3-geq-3}
\end{eqnarray}
It is clear to see that
%
\begin{eqnarray}
&&\sum _{k=3}^{\infty }f_{k+1}(2k+9)q^{k^{2}+9k+16}-\sum _{k=3}^{\infty }q^{k^{2}+7k+8}
\nonumber
\\[3pt]
&=&\sum _{k=4}^{\infty }f_{k}(2k+7)q^{k^{2}+7k+8}-\sum _{k=3}^{\infty }q^{k^{2}+7k+8}
\nonumber
\\[3pt]
&=&-q^{38}+\sum _{k=4}^{\infty }\left (f_{k}(2k+7)-1\right )q^{k^{2}+7k+8}.
\label{main-m-3-geq-3a}
\end{eqnarray}
Substituting  {\eqref{main-m-3-geq-3a}} into  {\eqref{main-m-3-geq-3}}, we obtain
\begin{eqnarray*}
\sum _{n= 0}^{\infty }S(n)q^{n}&=&-q^{38}+\sum _{k=3}^{\infty }q^{k^{2}+7k+7}
\left (1+\sum _{{n\geq 2 }\atop{n\neq 2k+9}}f_{k+1}(n)q^{n}\right )
\\[3pt]
&&+\sum _{k=4}^{\infty }\left (f_{k}(2k+7)-1\right )q^{k^{2}+7k+8}.
\end{eqnarray*}
From  {Theorem~\ref{lem-kernel-2}} (4), we see that $f_{k}(n)\geq 0$ for
$k\geq 4$ and $n\geq 2$ and $f_{k}(2k+7)\geq 1$. It follows that
$S(n)\geq 0$ when $n\geq 39$. Thus, by  {\eqref{main-m-3-claim-pf-re}}, we
conclude that $T_{2}(n)\geq 0$ for $n\geq 44$, and so
$M(1,n)-M(2,n)\geq 0$ for $n\geq 44$. Hence we complete the proof of  {Theorem~\ref{main-thm-n}} for $m=2$. \end{proof}

\section{On $M(m-1,n)\geq M(m,n)$ when $m\geq 3$}
\label{sec8}

In this section, we prove that  {Theorem~\ref{main-thm-n}} holds when
$m\geq 3$ by means of  {Theorem~\ref{lem-kernel}},  {Corollary~\ref{corollary}} and {Theorem~\ref{lem-4}}.

\begin{proof}[Proof of  {Theorem~\ref{main-thm-n}} for $m\geq 3$] Define
%
\begin{align}
\label{eq-def-tm}
\sum _{n=0}^{\infty }T_{m}(n)q^{n} :=&-q^{2m}+q^{2m+1}+q^{3m+1}+q^{m-1}
\frac{1-q}{(q^{2};q)_{m-2}}
\nonumber
\\[3pt]
&+\sum _{k=2}^{\infty
}\frac{q^{k(k+m)+3k+2m-2}}{(q^{3};q)_{k-2}(q^{2};q)_{k+m-2}}+\sum _{k=1}^{\infty
}\frac{q^{k(k+m)+4k+2m+2}(1-q^{m-2})}{(q^{2};q)_{k}(q^{2};q)_{k+m-2}} .\quad
\end{align}
From  {Theorem~\ref{lem-4}}, it can be seen that for $m\geq 3$ and
$n\geq 0$,
%
\begin{equation}
\label{mm-1-tm}
M(m-1,n)-M(m,n)\geq T_{m}(n).
\end{equation}
Define
%
\begin{equation}
\label{equ-def-um}
\sum _{n=0}^{\infty }U_{m}(n)q^{n}:= -q^{2m}+q^{2m+1}+q^{3m+1}+q^{m-1}
\frac{1-q}{(q^{2};q)_{m-2}}+\frac{q^{4m+8}}{(q^{2};q)_{m}},
\end{equation}
so
%
\begin{equation}
\label{tm-um}
\sum _{n=0}^{\infty }(T_{m}(n)-U_{m}(n))q^{n}=\sum _{k=3}^{\infty
}\frac{q^{k(k+m)+3k+2m-2}}{(q^{3};q)_{k-2}(q^{2};q)_{k+m-2}}+\sum _{k=1}^{\infty
}\frac{q^{k(k+m)+4k+2m+2}(1-q^{m-2})}{(q^{2};q)_{k}(q^{2};q)_{k+m-2}}.
\end{equation}

When $m=3$, observe that
\begin{equation*}
\sum _{k=1}^{\infty
}\frac{q^{k(k+m)+4k+2m+2}(1-q^{m-2})}{(q^{2};q)_{k}(q^{2};q)_{k+m-2}}=
\sum _{k=1}^{\infty
}\frac{q^{k^{2}+7k+8}(1-q)}{(q^{2};q)_{k}(q^{2};q)_{k+1}}=\sum _{n=0}^{\infty }T_{2}(n)q^{n+1},
\end{equation*}
where $T_{2}(n)$ is defined in  {\eqref{main-m-3-claim}}. From the proof of
 {Theorem~\ref{main-thm-n}} for $m=2$, we see that $T_{2}(n)\geq 0$ for
$n\geq 44$. Moreover, it can be checked that $T_{2}(43)\geq 0$.

When $m\geq 4$,
\begin{equation*}
\sum _{k=1}^{\infty
}\frac{q^{k(k+m)+4k+2m+2}(1-q^{m-2})}{(q^{2};q)_{k}(q^{2};q)_{k+m-2}} =
\sum _{k=1}^{\infty
}\frac{q^{k(k+m)+4k+2m+2}}{(q^{2};q)_{k}(q^{2};q)_{m-4}(q^{m-1};q)_{k+1}},
\end{equation*}
which obviously has nonnegative power series coefficients. So in either
case  {\eqref{tm-um}} implies that for $m\geq 3$ and $n\geq 44$,
%
\begin{equation}
\label{mm-1-tm-123}
T_{m}(n)\geq U_{m}(n).
\end{equation}

We proceed to establish the nonnegativity of $U_{m}(n)$ when
$m\geq 3$. Observe that
%
\begin{align}
\label{pf-main-case-5-99}
\sum _{n=0}^{\infty }U_{m}(n)q^{n} =&-q^{2m}+q^{2m+1}+q^{3m+1}+q^{m-1}
\left (1-q+\sum _{n=2}^{\infty }d_{m-1}(n)q^{n}\right )
\nonumber
\\[3pt]
&+\sum _{n=0}^{\infty }p_{m+1}(n)q^{n+4m+8}
\nonumber
\\[3pt]
=&q^{m-1}-q^{m}+q^{2m+1}+q^{3m+1}+\sum _{{n\geq 2}\atop{n\neq m+1}}d_{m-1}(n)q^{n+m-1}
\nonumber
\\[3pt]
&+(d_{m-1}(m+1)-1)q^{2m}+\sum _{n=0}^{\infty }p_{m+1}(n)q^{n+4m+8}.
\end{align}

There are two cases:

Case 1. When $3\leq m\leq 7$, from  {\eqref{pf-main-case-5-99}}, we see that for
$n\geq 44$,
%
\begin{equation}
\label{eq-um-d-p}
U_{m}(n)=d_{m-1}(n-m+1)+p_{m+1}(n-4m-8).
\end{equation}
By  {Theorem~\ref{lem-kernel}} (2)--(6), we see that when
$3\leq m\leq 7$ and $n\geq m-1$,
\begin{equation*}
d_{m-1}(n-m+1)\geq -\left \lfloor \frac{n-m+12}{12}\right \rfloor .
\end{equation*}
By  {Corollary~\ref{corollary}}, we have for $n\geq 4m+8$,
\begin{equation*}
p_{m+1}(n-4m-8)\geq \left \lfloor \frac{n-4m-8}{6}\right \rfloor .
\end{equation*}
Thus by  {\eqref{eq-um-d-p}}, we derive that when $3\leq m\leq 7$ and
$n\geq 72$,
\begin{equation*}
U_{m}(n)\geq \left \lfloor \frac{n-4m-8}{6}\right \rfloor -\left
\lfloor \frac{n-m+12}{12}\right \rfloor \geq 0.
\end{equation*}
It is trivial to check that $U_{m}(n)\geq 0$ for $44\leq n\leq 71$. So
we are led to $U_{m}(n)\geq 0$ for $n\geq 44$ and $3\leq m\leq 7$.

Case 2. When $m\geq 8$, from  {\eqref{pf-main-case-5-99}}, we derive that when
$n\geq m+1$ and $n\neq 2m$,
%
\begin{equation}
\label{eq-um-6-1}
U_{m}(n)\geq d_{m-1}(n-m+1)+p_{m+1}(n-4m-8).
\end{equation}
When $n=2m$,
%
\begin{equation}
\label{eq-um-6-2}
U_{m}(n)= d_{m-1}(m+1)-1.
\end{equation}
By  {Theorem~\ref{lem-kernel}} (7), we find that
$d_{m-1}(n-m+1)+p_{m+1}(n-4m-8)\geq 0$ and $d_{m-1}(m+1)\geq 1$ when
$m\geq 8$ and $n\geq m+1$. Hence by  {\eqref{eq-um-6-1}} and  {\eqref{eq-um-6-2}}, we derive that
$U_{m}(n)\geq 0$  when $m\geq 8$ and $n\geq m+1$. So $ U_{m}(n)\geq 0$  when $m\geq 3$ and
$n\geq 44$. Hence it follows from  {\eqref{mm-1-tm}} and  {\eqref{mm-1-tm-123}} that $M(m-1,n)-M(m,n)\geq 0$ for $m\geq 3$ and
$n\geq 44$. This completes the proof of  {Theorem~\ref{main-thm-n}} when
$m\geq 3$.\end{proof}

\section{On $M(0,n)\geq M(1,n)$}
\label{sec9}

In the section, we finish the proof of  {Theorem~\ref{main-thm-n}} by showing
that $M(0,n)\ge M(1,n)$ for $n\geq 44$. As has already been mentioned in Section~\ref{sec2}, the
proof of $M(0,n)\ge M(1,n)$ for $n\geq 44$ is the most complicated. Setting
$m=1$ in  {Theorem~\ref{lem-4-1}}, we find that
%
\begin{eqnarray}
&&\sum _{n=0}^{\infty }\left (M(0,n)-M(1,n)\right )q^{n}
\nonumber
\\[3pt]
&=&1-2q+q^{3}+q^{4}-q^{7}+\frac{q^{10}}{1-q^{2}}-
\frac{q^{11}}{1-q^{2}}
\nonumber
\\[3pt]
&& +\sum _{k=3}^{\infty
}\frac{q^{k^{2}+3k}}{(q^{2};q)_{k-1}(q^{2};q)_{k-3}}+\sum _{k=2}^{\infty }\frac{q^{k^{2}+4k}}{(q^{3};q)_{k-2}(q^{2};q)_{k-1}}
\nonumber
\\[3pt]
&&-\sum _{k=1}^{\infty
}\frac{q^{k^{2}+5k+3}}{(q^{2};q)_{k-1}(q^{2};q)_{k-1}}.
\label{lem-0-1-temp}
\end{eqnarray}

To show that  {Theorem~\ref{main-thm-n}} holds when $m=1$, we next aim to
show that the following generating function of $M(0,n)-M(1,n)$ holds.

\begin{thm}%
\label{lem-0-1}
We have
%
\begin{eqnarray}
&&\sum _{n=0}^{\infty }\left (M(0,n)-M(1,n)\right )q^{n}
\label{lem-0-1-e}
\\[3pt]
&=&1-2q+q^{3}+q^{4}-q^{7}-q^{9}+q^{18}+\frac{q^{10}}{1-q^{2}} -
\frac{q^{11}}{1-q^{2}}
\nonumber
\\[3pt]
&&+\frac{q^{12}}{1-q^{2}}-\frac{q^{17}}{1-q^{2}} -
\frac{q^{19}}{1-q^{2}}+\frac{q^{20}}{1-q^{2}} +\frac{q^{21}}{1-q^{2}}-
\frac{q^{23}}{(1-q^{2})^{2}}
\nonumber
\\[3pt]
&& +\frac{q^{24}}{(q^{2};q)_{2}}+\frac{q^{28}}{(1-q^{2})^{2}} -
\frac{q^{38}}{(1-q^{3})^{2}}+\sum _{k=3}^{\infty
}\frac{q^{k^{2}+5k}}{(q^{4};q)_{k-3}(q^{2};q)_{k-1}}
\nonumber
\\[3pt]
&&+\sum _{k=3}^{\infty
}\frac{q^{k^{2}+5k}}{(q^{3};q)_{k-2}(q^{2};q)_{k-3}}+\sum _{k= 3}^{\infty }\frac{q^{k^{2}+5k+2}}{(1-q^{2})(q^{4};q)_{k-3}(q^{2};q)_{k-3}}
\nonumber
\\[3pt]
&&+\sum _{k=3}^{\infty }\sum _{i=0}^{\infty
}\frac{q^{k^{2}+6k+5+(k-1)i}(1- q^{i+2})}{(q^{2};q)_{k-1}(q^{2};q)_{k-3}}+
\sum _{k= 3}^{\infty }\sum _{i= 0}^{\infty
}\frac{q^{k^{2}+9k+8+ik}(1-q^{i+8})}{(q^{2};q)_{k-2}(q^{3};q)_{k-2}(1-q^{k+1})}.
\nonumber
\end{eqnarray}
\end{thm}

 To obtain  {Theorem~\ref{lem-0-1}}, we are required to further
expand three summations in  {\eqref{lem-0-1-temp}}.

\begin{lem}%
\label{main-m-0-po-2-lem}%
We have
%
\begin{eqnarray}
&&\sum _{k=3}^{\infty
}\frac{q^{k^{2}+3k}}{(q^{2};q)_{k-1}(q^{2};q)_{k-3}}
\nonumber
\\[3pt]
&=&\sum _{k=3}^{\infty
}\frac{q^{k^{2}+3k}}{(q^{2};q)_{k-3}(q^{2};q)_{k-3}}+\sum _{k=3}^{\infty }\frac{q^{k^{2}+4k-1}}{(q^{2};q)_{k-2}(q^{2};q)_{k-3}}+\sum _{k=3}^{\infty }\frac{q^{k^{2}+4k}}{(q^{2};q)_{k-2}(q^{2};q)_{k-3}}
\nonumber
\\[3pt]
&&+\sum _{k=3}^{\infty
}\frac{q^{k^{2}+5k}}{(q^{2};q)_{k-1}(q^{2};q)_{k-3}}.
\label{main-m-0-po-2}
\end{eqnarray}
\end{lem}

\begin{proof} It is clear that
%
\begin{equation}
\label{lem-9-1}
\sum _{k=3}^{\infty
}\frac{q^{k^{2}+3k}}{(q^{2};q)_{k-1}(q^{2};q)_{k-3}}=\sum _{k=3}^{\infty }\frac{q^{k^{2}+3k}}{(q^{2};q)_{k-2}(q^{2};q)_{k-3}}\cdot
\frac{1}{1-q^{k}}.
\end{equation}
Obviously, when $k\geq 1$,
%
\begin{equation}
\label{equ-iden-1-1-qk}
\frac{1}{1-q^{k}}=(1-q^{k-1})+q^{k-1}+q^{k}+\frac{q^{2k}}{1-q^{k}},
\end{equation}
so we can obtain  {\eqref{main-m-0-po-2}} by substituting  {\eqref{equ-iden-1-1-qk}} into {\eqref{lem-9-1}}. \end{proof}

\begin{lem}
\label{main-m-0-po-1-lem}
We have
%
\begin{eqnarray}
&&\sum _{k=2}^{\infty
}\frac{q^{k^{2}+4k}}{(q^{3};q)_{k-2}(q^{2};q)_{k-1}}
\nonumber
\\[3pt]
=&&\frac{q^{12}}{1-q^{2}}+\sum _{k=3}^{\infty
}\frac{q^{k^{2}+4k}}{(q^{3};q)_{k-3}(q^{2};q)_{k-2}}+\sum _{k=3}^{\infty }\frac{q^{k^{2}+5k}}{(q^{3};q)_{k-2}(q^{2};q)_{k-1}}
\nonumber
\\[3pt]
&&+\sum _{k=3}^{\infty
}\frac{q^{k^{2}+5k}}{(q^{3};q)_{k-3}(q^{2};q)_{k-1}} .
\label{main-m-0-po-1}
\end{eqnarray}
\end{lem}

\begin{proof} Observe that
%
\begin{equation}
\label{eq-sum-k=2-inf-fac}
\sum _{k=2}^{\infty
}\frac{q^{k^{2}+4k}}{(q^{3};q)_{k-2}(q^{2};q)_{k-1}} =
\frac{q^{12}}{1-q^{2}}+\sum _{k=3}^{\infty
}\frac{q^{k^{2}+4k}}{(q^{3};q)_{k-3}(q^{2};q)_{k-1}}\cdot
\frac{1}{1-q^{k}}.
\end{equation}
Clearly, when $k\geq 1$,
%
\begin{equation}
\label{equ-fra-1-1-qk}
\frac{1}{1-q^{k}}=(1-q^{k})+\frac{q^{k}}{1-q^{k}}+q^{k}.
\end{equation}
Substituting  {\eqref{equ-fra-1-1-qk}} into  {\eqref{eq-sum-k=2-inf-fac}}, we
obtain  {\eqref{main-m-0-po-1}}. \end{proof}

\begin{lem}%
\label{lem-0-last}
We have
%
\begin{eqnarray}%
\label{equ-m-0-sum-3}
&&\sum _{k=1}^{\infty
}\frac{q^{k^{2}+5k+3}}{(q^{2};q)_{k-1}(q^{2};q)_{k-1}}
\nonumber
\\[3pt]
&=&q^{9}+\frac{q^{19}}{1-q^{2}}+\frac{q^{23}}{(1-q^{2})^{2}}+\sum _{k=
2}^{\infty }\frac{q^{k^{2}+5k+3}}{(q^{2};q)_{k-1}(q^{3};q)_{k-2}}+\sum _{k=
3}^{\infty }\frac{q^{k^{2}+5k+5}}{(q^{2};q)_{k-1}(q^{2};q)_{k-3}}
\nonumber
\\[3pt]
&&+\!\sum _{k=3}^{\infty
}\frac{q^{k^{2}+6k+4}}{(q^{2};q)_{k-1}(q^{2};q)_{k-2}}
+\!\sum _{k= 2}^{\infty }\frac{q^{k^{2}+6k+5}}{(q^{2};q)_{k-1}(q^{3};q)_{k-2}}+\!\sum _{k=
3}^{\infty
}\frac{q^{k^{2}+6k+7}}{(q^{2};q)_{k-1}(q^{2};q)_{k-3}(1-q^{k})}
\nonumber
\\[3pt]
&&+\sum _{k=3}^{\infty
}\frac{q^{k^{2}+7k+6}(1+q^{2})}{(q^{2};q)_{k-1}(q^{3};q)_{k-2}}+\sum _{k=
3}^{\infty }\frac{q^{k^{2}+7k+10}}{(q^{2};q)_{k-1}(q^{2};q)_{k-1}}.
\end{eqnarray}
\end{lem}

\begin{proof} Clearly,
%
\begin{equation}
\label{equ-sum-k=1-inf-fac}
\sum _{k=1}^{\infty
}\frac{q^{k^{2}+5k+3}}{(q^{2};q)_{k-1}(q^{2};q)_{k-1}}=q^{9}+\sum _{k=2}^{\infty }\frac{q^{k^{2}+5k+3}}{(q^{2};q)_{k-1}(q^{3};q)_{k-2}}\cdot
\frac{1}{1-q^{2}}.
\end{equation}
It is trivial to check that
%
\begin{equation}
\label{equ-fac-1-1-q2-1+fac}
\frac{1}{1-q^{2}}=1+\frac{q^{2}(1-q^{k})}{1-q^{2}}+
\frac{q^{k+2}}{1-q^{2}}.
\end{equation}
Substituting  {\eqref{equ-fac-1-1-q2-1+fac}} into  {\eqref{equ-sum-k=1-inf-fac}}, we have
%
\begin{eqnarray}
\sum _{k=1}^{\infty
}\frac{q^{k^{2}+5k+3}}{(q^{2};q)_{k-1}(q^{2};q)_{k-1}} &=&q^{9}+\sum _{k=2}^{\infty }\frac{q^{k^{2}+5k+3}}{(q^{2};q)_{k-1}(q^{3};q)_{k-2}}+\sum _{k=2}^{\infty }\frac{q^{k^{2}+5k+5}}{(q^{2};q)_{k-1}(q^{2};q)_{k-2}}
\nonumber
\\[3pt]
&&+\sum _{k=2}^{\infty
}\frac{q^{k^{2}+6k+5}}{(q^{2};q)_{k-1}(q^{2};q)_{k-1}}.
\label{equ-q9+sum-k-2-inf}
\end{eqnarray}
We next transform the third and the fourth terms in  {\eqref{equ-q9+sum-k-2-inf}} respectively. First, we rewrite the
third term as follows.
%
\begin{eqnarray}
&&\sum _{k=2}^{\infty
}\frac{q^{k^{2}+5k+5}}{(q^{2};q)_{k-1}(q^{2};q)_{k-2}}
\nonumber
\\[3pt]
&=&\frac{q^{19}}{1-q^{2}}+ \sum _{k= 3}^{\infty
}\frac{q^{k^{2}+5k+5}}{(q^{2};q)_{k-1}(q^{2};q)_{k-3}}\left (1+
\frac{q^{k-1}}{1-q^{k-1}} \right )
\nonumber
\\[3pt]
&=&\frac{q^{19}}{1-q^{2}}+ \sum _{k=3}^{\infty
}\frac{q^{k^{2}+5k+5}}{(q^{2};q)_{k-1}(q^{2};q)_{k-3}}+\sum _{k=3}^{\infty }\frac{q^{k^{2}+6k+4}}{(q^{2};q)_{k-1}(q^{2};q)_{k-2}}.
\label{equ-su-k-2-inf}
\end{eqnarray}
We next transform the last term of  {\eqref{equ-q9+sum-k-2-inf}} as given below.
%
\begin{align}
&\sum _{k=2}^{\infty
}\frac{q^{k^{2}+6k+5}}{(q^{2};q)_{k-1}(q^{2};q)_{k-1}}
\nonumber
\\[3pt]
=&\sum _{k= 2}^{\infty
}\frac{q^{k^{2}+6k+5}}{(q^{2};q)_{k-1}(q^{3};q)_{k-2}}\left (1+
\frac{q^{2}}{1-q^{2}}\right )
\nonumber
\\[3pt]
=&\sum _{k=2}^{\infty
}\frac{q^{k^{2}+6k+5}}{(q^{2};q)_{k-1}(q^{3};q)_{k-2}}+\sum _{k= 2}^{\infty }\frac{q^{k^{2}+6k+7}}{(q^{2};q)_{k-1}(q^{2};q)_{k-1}}
\nonumber
\\[3pt]
=&\sum _{k= 2}^{\infty
}\frac{q^{k^{2}+6k+5}}{(q^{2};q)_{k-1}(q^{3};q)_{k-2}}+
\frac{q^{23}}{(1-q^{2})^{2}}+\sum _{k= 3}^{\infty
}\frac{q^{k^{2}+6k+7}}{(q^{2};q)_{k-1}(q^{2};q)_{k-3}(1-q^{k})}\cdot
\frac{1}{1-q^{k-1}}.
\label{equ-su-k-2-in-far=q=k2}
\end{align}
Note that
%
\begin{equation}
\label{eq-fac-1-1-qk-1}
\frac{1}{1-q^{k-1}}=1+\frac{q^{k-1}(1-q^{4})}{1-q^{k-1}}+
\frac{q^{k+3}}{1-q^{k-1}}.
\end{equation}
Substituting  {\eqref{eq-fac-1-1-qk-1}} into  {\eqref{equ-su-k-2-in-far=q=k2}} and with some simplification, we deduce
that
%
\begin{eqnarray}
&&\sum _{k=2}^{\infty
}\frac{q^{k^{2}+6k+5}}{(q^{2};q)_{k-1}(q^{2};q)_{k-1}}
\nonumber
\\[3pt]
&=&\frac{q^{23}}{(1-q^{2})^{2}}+\sum _{k=2}^{\infty
}\frac{q^{k^{2}+6k+5}}{(q^{2};q)_{k-1}(q^{3};q)_{k-2}}+\sum _{k= 3}^{\infty }\frac{q^{k^{2}+6k+7}}{(q^{2};q)_{k-1}(q^{2};q)_{k-3}(1-q^{k})}
\nonumber
\\[3pt]
&&+\sum _{k=3}^{\infty
}\frac{q^{k^{2}+7k+6}(1+q^{2})}{(q^{2};q)_{k-1}(q^{3};q)_{k-2}}+\sum _{k=
3}^{\infty }\frac{q^{k^{2}+7k+10}}{(q^{2};q)_{k-1}(q^{2};q)_{k-1}}.
\label{main-m-0-ne-case-5-6}
\end{eqnarray}
Substituting  {\eqref{equ-su-k-2-inf}} and  {\eqref{main-m-0-ne-case-5-6}} into {\eqref{equ-q9+sum-k-2-inf}}, we obtain  {\eqref{equ-m-0-sum-3}}. Thus the lemma
has been verified.\end{proof}

We are now in a position to give a proof of  {Theorem~\ref{lem-0-1}} in light
of  {Lemmas~\ref{main-m-0-po-2-lem}, \ref{main-m-0-po-1-lem} and \ref{lem-0-last}}.

\begin{proof}[Proof of  {Theorem~\ref{lem-0-1}}] Substituting  {\eqref{main-m-0-po-2}}, {\eqref{main-m-0-po-1}} and  {\eqref{equ-m-0-sum-3}} into  {\eqref{lem-0-1-temp}}, we arrive at
%
\begin{eqnarray}
&&\sum _{n=0}^{\infty }\left (M(0,n)-M(1,n)\right )q^{n}
\label{equ-gf-m1-case-final}
\\[3pt]
&=&1-2q+q^{3}+q^{4}-q^{7}-q^{9}+\frac{q^{10}}{1-q^{2}}-
\frac{q^{11}}{1-q^{2}}+\frac{q^{12}}{1-q^{2}}-\frac{q^{19}}{1-q^{2}}-
\frac{q^{23}}{(1-q^{2})^{2}}
\nonumber
\\[3pt]
&&+\left (\sum _{k=3}^{\infty
}\frac{q^{k^{2}+4k}}{(q^{2};q)_{k-2}(q^{3};q)_{k-3}}-\sum _{k=2}^{\infty }\frac{q^{k^{2}+6k+5}}{(q^{2};q)_{k-1}(q^{3};q)_{k-2}}\right )
\label{gf-m-1-case-6}
\\[3pt]
&&+\left (\sum _{k=3}^{\infty
}\frac{q^{k^{2}+4k-1}}{(q^{2};q)_{k-2}(q^{2};q)_{k-3}}-\sum _{k=3}^{\infty }\frac{q^{k^{2}+6k+4}}{(q^{2};q)_{k-1}(q^{2};q)_{k-2}}\right )
\label{gf-m-1-case-5}
\\[3pt]
&&+\left (\sum _{k=3}^{\infty
}\frac{q^{k^{2}+3k}}{(q^{2};q)_{k-3}(q^{2};q)_{k-3}}-\sum _{k=3}^{\infty }\frac{q^{k^{2}+7k+10}}{(q^{2};q)_{k-1}(q^{2};q)_{k-1}}\right )
\label{gf-m-1-case-7}
\\[3pt]
&&+\left (\sum _{k=3}^{\infty
}\frac{q^{k^{2}+5k}}{(q^{2};q)_{k-1}(q^{3};q)_{k-2}}-\sum _{k=2}^{\infty }\frac{q^{k^{2}+5k+3}}{(q^{2};q)_{k-1}(q^{3};q)_{k-2}}\right )
\label{gf-m-1-case-3}
\\[3pt]
&&+\left (\sum _{k=3}^{\infty
}\frac{q^{k^{2}+5k}}{(q^{2};q)_{k-1}(q^{2};q)_{k-3}}-\sum _{k=3}^{\infty }\frac{q^{k^{2}+5k+5}}{(q^{2};q)_{k-1}(q^{2};q)_{k-3}}\right )
\label{gf-m-1-case-4}
\\[3pt]
&&+\left (\sum _{k=3}^{\infty
}\frac{q^{k^{2}+4k}}{(q^{2};q)_{k-2}(q^{2};q)_{k-3}}-\sum _{k=3}^{\infty }\frac{q^{k^{2}+6k+7}}{(q^{2};q)_{k-1}(q^{2};q)_{k-3}(1-q^{k})}
\right )
\label{gf-m-1-case-8}
\\[3pt]
&&+\left (\sum _{k=3}^{\infty
}\frac{q^{k^{2}+5k}}{(q^{2};q)_{k-1}(q^{3};q)_{k-3}}-\sum _{k=3}^{\infty }\frac{q^{k^{2}+7k+6}(1+q^{2})}{(q^{2};q)_{k-1}(q^{3};q)_{k-2}}
\right )
\label{gf-m-1-case-9}
.
\label{main-m-0-dis}
\end{eqnarray}
We proceed to simplify seven differences in the above identity. Note that
\begin{equation*}
\sum _{k=3}^{\infty
}\frac{q^{k^{2}+4k}}{(q^{2};q)_{k-2}(q^{3};q)_{k-3}}=\sum _{k=2}^{\infty }\frac{q^{k^{2}+6k+5}}{(q^{2};q)_{k-1}(q^{3};q)_{k-2}},
\end{equation*}
so  {\eqref{gf-m-1-case-6}} is equal to $0$.

We now consider the difference  {\eqref{gf-m-1-case-5}}. Observe that
\begin{equation*}
\sum _{k=3}^{\infty
}\frac{q^{k^{2}+4k-1}}{(q^{2};q)_{k-2}(q^{2};q)_{k-3}}=\sum _{k=2}^{\infty }\frac{q^{k^{2}+6k+4}}{(q^{2};q)_{k-1}(q^{2};q)_{k-2}}.
\end{equation*}
Hence  {\eqref{gf-m-1-case-5}} is equal to
%
\begin{equation}
\label{main-m-0-5}
\sum _{k=2}^{\infty
}\frac{q^{k^{2}+6k+4}}{(q^{2};q)_{k-1}(q^{2};q)_{k-2}}-\sum _{k=3}^{\infty }\frac{q^{k^{2}+6k+4}}{(q^{2};q)_{k-1}(q^{2};q)_{k-2}} =
\frac{q^{20}}{1-q^{2}}.
\end{equation}
Note that
\begin{equation*}
\sum _{k=3}^{\infty
}\frac{q^{k^{2}+3k}}{(q^{2};q)_{k-3}(q^{2};q)_{k-3}}=\sum _{k=1}^{\infty }\frac{q^{k^{2}+7k+10}}{(q^{2};q)_{k-1}(q^{2};q)_{k-1}},
\end{equation*}
so  {\eqref{gf-m-1-case-7}} is equal to
%
\begin{equation}
\label{main-m-0-7}
\sum _{k=1}^{\infty
}\frac{q^{k^{2}+7k+10}}{(q^{2};q)_{k-1}(q^{2};q)_{k-1}}-\sum _{k=3}^{\infty }\frac{q^{k^{2}+7k+10}}{(q^{2};q)_{k-1}(q^{2};q)_{k-1}} =q^{18}+
\frac{q^{28}}{(1-q^{2})^{2}}.
\end{equation}

For  {\eqref{gf-m-1-case-3}}, we see that
%
\begin{eqnarray}
&&\sum _{k=3}^{\infty
}\frac{q^{k^{2}+5k}}{(q^{2};q)_{k-1}(q^{3};q)_{k-2}}-\sum _{k=2}^{\infty }\frac{q^{k^{2}+5k+3}}{(q^{2};q)_{k-1}(q^{3};q)_{k-2}}
\nonumber
\\[3pt]
&=&\sum _{k=3}^{\infty
}\frac{q^{k^{2}+5k}}{(q^{2};q)_{k-1}(q^{3};q)_{k-2}}-
\frac{q^{17}}{1-q^{2}}-\sum _{k=3}^{\infty
}\frac{q^{k^{2}+5k+3}}{(q^{2};q)_{k-1}(q^{3};q)_{k-2}}
\nonumber
\\[3pt]
&=&-\frac{q^{17}}{1-q^{2}}+\sum _{k=3}^{\infty
}\frac{q^{k^{2}+5k}(1-q^{3})}{(q^{2};q)_{k-1}(q^{3};q)_{k-2}}
\nonumber
\\[3pt]
&=&-\frac{q^{17}}{1-q^{2}}+\sum _{k=3}^{\infty
}\frac{q^{k^{2}+5k}}{(q^{4};q)_{k-3}(q^{2};q)_{k-1}}.
\label{main-m-0-3}
\end{eqnarray}
The difference  {\eqref{gf-m-1-case-4}} can be simplified as follows:
%
\begin{eqnarray}
&&\sum _{k=3}^{\infty
}\frac{q^{k^{2}+5k}}{(q^{2};q)_{k-1}(q^{2};q)_{k-3}}-\sum _{k=3}^{\infty }\frac{q^{k^{2}+5k+5}}{(q^{2};q)_{k-1}(q^{2};q)_{k-3}}
\nonumber
\\[3pt]
&=& {\sum _{k=3}^{\infty
}\frac{q^{k^{2}+5k}(1-q^{5})}{(q^{2};q)_{k-1}(q^{2};q)_{k-3}}}
\nonumber
\\[3pt]
&=&\sum _{k=3}^{\infty
}\frac{q^{k^{2}+5k}(1-q^{2}+q^{2}-q^{5})}{(q^{2};q)_{k-1}(q^{2};q)_{k-3}}
\nonumber
\\[3pt]
&=&\sum _{k=3}^{\infty
}\frac{q^{k^{2}+5k}(1-q^{2})}{(q^{2};q)_{k-1}(q^{2};q)_{k-3}}+\sum _{k=3}^{\infty }\frac{q^{k^{2}+5k+2}(1-q^{3})}{(q^{2};q)_{k-1}(q^{2};q)_{k-3}}
\nonumber
\\[3pt]
&=&\sum _{k=3}^{\infty
}\frac{q^{k^{2}+5k}}{(q^{3};q)_{k-2}(q^{2};q)_{k-3}}+\sum _{k=3}^{\infty }\frac{q^{k^{2}+5k+2}}{(1-q^{2})(q^{4};q)_{k-3}(q^{2};q)_{k-3}}.
\label{main-m-0-4}
\end{eqnarray}
For  {\eqref{gf-m-1-case-8}}, we have
%
\begin{eqnarray}%
\label{main-m-0-8}
&&\sum _{k=3}^{\infty
}\frac{q^{k^{2}+4k}}{(q^{2};q)_{k-2}(q^{2};q)_{k-3}}-\sum _{k=3}^{\infty }\frac{q^{k^{2}+6k+7}}{(q^{2};q)_{k-1}(q^{2};q)_{k-3}(1-q^{k})}
\nonumber
\\[3pt]
&=& {\sum _{k=2}^{\infty
}\frac{q^{k^{2}+6k+5}}{(q^{2};q)_{k-1}(q^{2};q)_{k-2}}-\sum _{k=3}^{\infty }\frac{q^{k^{2}+6k+7}}{(q^{2};q)_{k-1}(q^{2};q)_{k-3}(1-q^{k})}}
\nonumber
\\[3pt]
&=&\frac{q^{21}}{1-q^{2}}+\sum _{k=3}^{\infty
}\frac{q^{k^{2}+6k+5}}{(q^{2};q)_{k-1}(q^{2};q)_{k-2}}-\sum _{k=3}^{\infty }\frac{q^{k^{2}+6k+7}}{(q^{2};q)_{k-1}(q^{2};q)_{k-3}(1-q^{k})}
\nonumber
\\[3pt]
&=&\frac{q^{21}}{1-q^{2}}+\sum _{k=3}^{\infty
}\frac{q^{k^{2}+6k+5}}{(q^{2};q)_{k-1}(q^{2};q)_{k-3}}\left (
\frac{1}{1-q^{k-1}} -\frac{q^{2}}{1-q^{k}}\right )
\nonumber
\\[3pt]
&=&\frac{q^{21}}{1-q^{2}}+\sum _{k=3}^{\infty
}\frac{q^{k^{2}+6k+5}}{(q^{2};q)_{k-1}(q^{2};q)_{k-3}}\left (\sum _{i=0}^{\infty }q^{(k-1)i} - {\sum _{i=0}^{\infty }q^{ki+2}}\right )
\nonumber
\\[3pt]
&=&\frac{q^{21}}{1-q^{2}}+\sum _{k=3}^{\infty }\sum _{i=0}^{\infty
}\frac{q^{k^{2}+6k+5+(k-1)i}(1- q^{i+2})}{(q^{2};q)_{k-1}(q^{2};q)_{k-3}}.
\end{eqnarray}

Finally, we transform  {\eqref{gf-m-1-case-9}} as given below:
%
\begin{eqnarray}%
\label{eq-sum=k=3-in-fac-q-k2}
&&\sum _{k=3}^{\infty
}\frac{q^{k^{2}+5k}}{(q^{2};q)_{k-1}(q^{3};q)_{k-3}}-\sum _{k=3}^{\infty }\frac{q^{k^{2}+7k+6}(1+q^{2})}{(q^{2};q)_{k-1}(q^{3};q)_{k-2}}
\nonumber
\\[3pt]
&=&\left (\frac{q^{24}}{(q^{2};q)_{2}}+\sum _{k=3}^{\infty
}\frac{q^{k^{2}+7k+6}}{(q^{2};q)_{k}(q^{3};q)_{k-2}}\right )-\sum _{k=3}^{\infty }\frac{q^{k^{2}+7k+6}(1+q^{2})}{(q^{2};q)_{k-1}(q^{3};q)_{k-2}}
\nonumber
\\[3pt]
&=&\frac{q^{24}}{(q^{2};q)_{2}}+\sum _{k=3}^{\infty
}\frac{q^{k^{2}+7k+6}}{(q^{2};q)_{k-1}(q^{3};q)_{k-2}}\left (
\frac{1}{1-q^{k+1}}-(1+q^{2})\right )
\nonumber
\\[3pt]
&=&\frac{q^{24}}{(q^{2};q)_{2}}+\sum _{k=3}^{\infty
}\frac{q^{k^{2}+7k+6}}{(q^{2};q)_{k-1}(q^{3};q)_{k-2}}\left (
\frac{q^{2k+2}}{1-q^{k+1}}-q^{2}(1-q^{k-1})\right )
\nonumber
\\[3pt]
&=&\frac{q^{24}}{(q^{2};q)_{2}}+\sum _{k=3}^{\infty
}\frac{q^{k^{2}+9k+8}}{(q^{2};q)_{k}(q^{3};q)_{k-2}}-\sum _{k=3}^{\infty }\frac{q^{k^{2}+7k+8}}{(q^{2};q)_{k-3}(1-q^{k})(q^{3};q)_{k-2}}
\nonumber
\\[3pt]
&=&\frac{q^{24}}{(q^{2};q)_{2}}+\sum _{k=3}^{\infty
}\frac{q^{k^{2}+9k+8}}{(q^{2};q)_{k}(q^{3};q)_{k-2}}-\sum _{k=2}^{\infty
}\frac{q^{k^{2}+9k+16}}{(q^{2};q)_{k-2}(1-q^{k+1})(q^{3};q)_{k-1}}
\nonumber
\\[3pt]
&=&\frac{q^{24}}{(q^{2};q)_{2}}-\frac{q^{38}}{(1-q^{3})^{2}}+\sum _{k=3}^{\infty
}\frac{q^{k^{2}+9k+8}}{(q^{2};q)_{k-2}(q^{3};q)_{k-2}(1-q^{k+1})}
\left (\frac{1}{1-q^{k}}-\frac{q^{8}}{1-q^{k+1}}\right )
\nonumber
\\[3pt]
&=&\frac{q^{24}}{(q^{2};q)_{2}}-\frac{q^{38}}{(1-q^{3})^{2}}+\sum _{k=3}^{\infty
}\frac{q^{k^{2}+9k+8}}{(q^{2};q)_{k-2}(q^{3};q)_{k-2}(1-q^{k+1})}
\left (\sum _{i=0}^{\infty }q^{ik}-\sum _{i=0}^{\infty }q^{i(k+1)+8}
\right )
\nonumber
\\[3pt]
&=&\frac{q^{24}}{(q^{2};q)_{2}}-\frac{q^{38}}{(1-q^{3})^{2}}+\sum _{k=3}^{\infty }\sum _{i= 0}^{\infty
}\frac{q^{k^{2}+9k+8+ik}(1-q^{i+8})}{(q^{2};q)_{k-2}(q^{3};q)_{k-2}(1-q^{k+1})}.
\label{main-m-0-9}
\end{eqnarray}
Substituting  {\eqref{main-m-0-5}} $\sim $  {\eqref{main-m-0-9}} into {\eqref{equ-gf-m1-case-final}}, we obtain  {\eqref{lem-0-1-e}}. This completes
the proof. \end{proof}

We are now ready to show that  {Theorem~\ref{main-thm-n}} holds when
$m=1$ with the aid of the generating function of $M(0,n)-M(1,n)$ in  {Theorem~\ref{lem-0-1}}.

\begin{proof}[Proof of  {Theorem~\ref{main-thm-n}} for $m=1$] Define
%
\begin{eqnarray}%
\label{equ-t1n-q11-q17}
\sum _{n=0}^{\infty }T_{1}(n)q^{n}&:=&-\frac{q^{11}}{1-q^{2}} -
\frac{q^{17}}{1-q^{2}}-\frac{q^{19}}{1-q^{2}} -
\frac{q^{23}}{(1-q^{2})^{2}}-\frac{q^{38}}{(1-q^{3})^{2}}
\nonumber
\\[3pt]
&&+\sum _{k=3}^{\infty
}\frac{q^{k^{2}+5k}}{(q^{4};q)_{k-3}(q^{2};q)_{k-1}}+\sum _{k=3}^{\infty }\sum _{i=0}^{\infty
}\frac{q^{k^{2}+6k+5+(k-1)i}(1- q^{i+2})}{(q^{2};q)_{k-1}(q^{2};q)_{k-3}}
\nonumber
\\[3pt]
&&+\sum _{k=3}^{\infty }\sum _{i= 0}^{\infty
}\frac{q^{k^{2}+9k+8+ik}(1-q^{i+8})}{(q^{2};q)_{k-2}(q^{3};q)_{k-2}(1-q^{k+1})}.
\end{eqnarray}
By  {Theorem~\ref{lem-0-1}}, we see that when $n\geq {10}$,
%
\begin{equation}
\label{equ-0-gf-final-two-summ}
M(0,n)-M(1,n)\geq T_{1}(n).
\end{equation}
To prove  {Theorem~\ref{main-thm-n}} for $m=1$, it suffices to establish the nonnegativity on $T_{1}(n)$.

By  {Lemma~\ref{le-chanmao-10}}, we see that
\begin{equation*}
\sum _{k=3}^{\infty }\sum _{i=0}^{\infty
}\frac{q^{k^{2}+6k+5+(k-1)i}(1- q^{i+2})}{(q^{2};q)_{k-1}(q^{2};q)_{k-3}}=
\sum _{k=3}^{\infty }\sum _{i=0}^{\infty
}\frac{q^{k^{2}+6k+5+(k-1)i}}{(q^{4};q)_{k-3}(q^{2};q)_{k-3}}\cdot
\frac{ 1- q^{i+2}}{(1-q^{2})(1-q^{3})}
\end{equation*}
and
\begin{align*}
&\sum _{k=3}^{\infty }\sum _{i= 0}^{\infty
}\frac{q^{k^{2}+9k+8+ik}(1-q^{i+8})}{(q^{2};q)_{k-2}(q^{3};q)_{k-2}(1-q^{k+1})}
\\
& =\sum _{k=3}^{\infty }\sum _{i=0}^{\infty
}\frac{q^{k^{2}+9k+8+ik}}{(q^{3};q)_{k-3}(q^{4};q)_{k-3}(1-q^{k+1})}
\cdot \frac{1-q^{i+8} }{(1-q^{2})(1-q^{3})}
\end{align*}
have nonnegative power series coefficients. Define
%
\begin{equation}
\label{equ-0-gf-final-two-summ-1-tema}
\sum _{n= 0}^{\infty }H(n)q^{n}:=\frac{q^{36}}{(q^{2};q)_{3}}-
\frac{q^{11}}{1-q^{2}} -\frac{q^{17}}{1-q^{2}}-\frac{q^{19}}{1-q^{2}} -
\frac{q^{23}}{(1-q^{2})^{2}} -\frac{q^{38}}{(1-q^{3})^{2}}.
\end{equation}
Note that
\begin{eqnarray*}
&&\sum _{k=3}^{\infty
}\frac{q^{k^{2}+5k}}{(q^{4};q)_{k-3}(q^{2};q)_{k-1}}-
\frac{q^{36}}{(q^{2};q)_{3}}
\\[3pt]
&&=\frac{q^{24}}{(q^{2};q)_{2}}+\frac{q^{36}}{(1-q^{4})(q^{2};q)_{3}}-
\frac{q^{36}}{(q^{2};q)_{3}}+\sum _{k=5}^{\infty
}\frac{q^{k^{2}+5k}}{(q^{4};q)_{k-3}(q^{2};q)_{k-1}}
\\[3pt]
&&=\frac{q^{24}}{(q^{2};q)_{2}}+\frac{q^{40}}{(1-q^{4})(q^{2};q)_{3}}+
\sum _{k=5}^{\infty
}\frac{q^{k^{2}+5k}}{(q^{4};q)_{k-3}(q^{2};q)_{k-1}},
\end{eqnarray*}
which has nonnegative power series coefficients. Hence we derive that for
$n\geq 11$,
%
\begin{equation}
\label{equ-0-gf-final-two-summ-1-tem}
T_{1}(n)\geq H(n).
\end{equation}
We proceed to show that $H(n)\geq 0$ when $n\geq 106$. By  {\eqref{gf-pkn}}, we see that
\begin{equation*}
\frac{q^{36}}{(q^{2};q)_{3}}=\sum _{n=36}^{\infty }p_{4}(n-36)q^{n}.
\end{equation*}
From  {\eqref{equ-0-gf-final-two-summ-1-tema}}, we find that for
$n\geq 38$,
\begin{equation*}
H(n)=
\begin{cases}
p_{4}(n-36),&\text{if }n\equiv 0,4\pmod{6},
\\[3pt]
p_{4}(n-36)-\frac{n-35}{3},&\text{if }n\equiv 2\pmod{6},
\\[3pt]
p_{4}(n-36)-\frac{n-21}{2}-3,&\text{if }n\equiv 1,3\pmod{6},
\\[3pt]
p_{4}(n-36)-\frac{n-21}{2}-3-\frac{n-35}{3},&\text{if }n\equiv 5
\pmod{6}.
\end{cases}
\end{equation*}
By  {Lemma~\ref{lem-kernel-f}} (3), we see that for $n\geq 48$,
\begin{equation*}
p_{4}(n-36)\geq 3\left (\frac{n-36}{12}-1\right )^{2}=
\frac{n^{2}}{48}-2n+48.
\end{equation*}
Hence we deduce that for $n\geq 106$,
\begin{equation*}
H(n)\geq {p_{4}(n-36)-\frac{n-21}{2}-3-\frac{n-35}{3}\geq }
\frac{n^{2}}{48}-2n+48-\frac{n-21}{2}-3-\frac{n-35}{3}\geq 0.
\end{equation*}
This implies that $T_{1}(n)\geq 0$ for $n\geq 106$, and so
$M(0,n)-M(1,n)\geq 0$ when $n\geq 106$. It can be checked that
$M(0,n)\geq M(1,n)$ when $44\leq n\leq 105$. Thus, we complete the proof
of the theorem. \end{proof}

\section{Proofs of  {Theorem~\ref{lem-indu}} and  {Conjecture~\ref{conj-opst3}}}
\label{sec10}

In this section, we first prove  {Theorem~\ref{lem-indu}}, and then give a
proof of  {Conjecture~\ref{conj-opst3}} with the aid of  {Theorem~\ref{main-thm-n}} and {Theorem~\ref{lem-indu}}.

To prove  {Theorem~\ref{lem-indu}}, setting $m=0$ in  {Theorem~\ref{thm-lem-gf-c-1}}, we see that
%
\begin{equation}
\label{gf-pn-21mn0-1-21}
\sum _{n= 0}^{\infty }21M(0,n)q^{n}=21-21q+\sum _{k= 1}^{\infty
}\frac{21q^{k^{2}+2k}}{(q;q)_{k}(q^{2};q)_{k-1}}.
\end{equation}
From \cite[Corollary 2.6]{and76},
%
\begin{equation}
\label{gf-pn-21mn0-2}
\sum _{n=0}^{\infty }p(n)q^{n}=1+\sum _{k=1}^{\infty
}\frac{q^{k^{2}}}{(q;q)_{k}^{2}}.
\end{equation}
Subtracting  {\eqref{gf-pn-21mn0-1-21}} from  {\eqref{gf-pn-21mn0-2}}, we obtain
the following generating function:
\begin{equation*}
\sum _{n=0}^{\infty }(p(n)-21M(0,n))q^{n}=-20+21q+\sum _{k= 1}^{\infty }q^{k^{2}}
\left (\frac{1}{(q;q)_{k}^{2}}-
\frac{21q^{2k}}{(q;q)_{k}(q^{2};q)_{k-1}} \right ).
\end{equation*}
For $k\geq 1$, define
%
\begin{equation}
\label{def-gmn}
\sum _{n= 0}^{\infty }g_{k}(n)q^{n}:=\frac{1}{(q;q)_{k}^{2}},
\end{equation}
and
%
\begin{equation}
\label{def-hmn}
\sum _{n= 0}^{\infty }h_{k}(n)q^{n}:=
\frac{q^{2k}}{(q;q)_{k}(q^{2};q)_{k-1}}.
\end{equation}
This leads to for $n\geq 2$,
%
\begin{equation}
\label{re-lem-indu}
p(n)-21M(0,n)=\sum _{k=1}^{\infty }\left (g_{k}(n-k^{2})-21h_{k}(n-k^{2})
\right ).
\end{equation}

The following theorem establishes the nonnegativity of
$g_{k}(n)-21h_{k}(n)$ which implies that $p(n)\geq 21 M(0,n)$ for
$n\geq 76$. Furthermore, it is not difficult to check that
$p(n)\geq 21 M(0,n)$ for $39\leq n\leq 75$. Hence by  {\eqref{re-lem-indu}}, we see that {Theorem~\ref{lem-indu}} immediately follows from the following theorem.

\begin{thm}%
\label{lem-indu-im}
\begin{itemize}
\item[\textup{(1)}] $g_{1}(n)\geq 21 h_{1}(n)$ for $n\geq 20$.
\item[\textup{(2)}] $g_{2}(n)\geq 21 h_{2}(n)$ for $n\geq 51$.
\item[\textup{(3)}] $g_{3}(n)\geq 21 h_{3}(n)$ for $n\geq 67$.
\item[\textup{(4)}] When $k\geq 4$, $g_{k}(n)\geq 21h_{k}(n)$ for
$n\geq 0$.
\end{itemize}
\end{thm}

Before proving  {Theorem~\ref{lem-indu-im}}, we first derive the following
recurrences of $g_{k}(n)$ and $h_{k}(n)$.

\begin{lem}%
\label{lem-rec-fmn}
For $k\geq 1$,
%
\begin{equation}
\label{rec-fmn}
g_{k}(n)=\sum _{i=0}^{\left \lfloor \frac{n}{k}\right \rfloor }(i+1)g_{k-1}(n-ki),
\end{equation}
and for $k \geq 2$,
%
\begin{equation}
\label{rec-gmn}
h_{k}(n)=\sum _{i=0}^{\left \lfloor \frac{n}{k}\right \rfloor -2}(i+1)h_{k-1}(n-ki-2).
\end{equation}
\end{lem}

\begin{proof} From the definition  {\eqref{def-gmn}} of $g_{k}(n)$, we see
that when $k\geq 1$,
\begin{eqnarray*}
\sum _{n=0}^{\infty }g_{k}(n)q^{n}&=&\frac{1}{(q;q)_{k}^{2}}
\\[3pt]
&=&\frac{1}{(q;q)_{k-1}^{2}}\cdot \frac{1}{(1-q^{k})^{2}}
\\[3pt]
&=&\sum _{n=0}^{\infty }g_{k-1}(n)q^{n}\cdot \sum _{i=0}^{\infty }(i+1)q^{ki}
\\[3pt]
&=&\sum _{n=0}^{\infty }\sum _{i=0}^{\left \lfloor \frac{n}{k}\right
\rfloor }(i+1)g_{k-1}(n-ki)q^{n}.
\end{eqnarray*}
So we obtain the recurrence  {\eqref{rec-fmn}} by equating coefficients of
$q^{n}$ on both sides of the above identity.

Proceeding as in the proof of  {\eqref{rec-fmn}}, we have
%
\begin{equation}
\label{rec-gmn-tem}
h_{k}(n)=\sum _{i=0}^{\infty }(i+1)h_{k-1}(n-ki-2).
\end{equation}
From  {\eqref{def-hmn}}, we see that $h_{k}(n)=0$ if $n<2k$. Thus  {\eqref{rec-gmn-tem}} can be written as follows:
\begin{equation*}
h_{k}(n)=\sum _{i=0}^{\left \lfloor \frac{n}{k}\right \rfloor -2 }(i+1)h_{k-1}(n-ki-2),
\end{equation*}
which is  {\eqref{rec-gmn}}. This completes the proof.\end{proof}

In order to prove  {Theorem~\ref{lem-indu-im}}, we also require the following
lemma.

\begin{lem}
When $k\geq 1$ and $n\geq 0$,
%
\begin{equation}
\label{mono-gk}
g_{k}(n+1)\geq g_{k}(n),
\end{equation}
and
%
\begin{equation}
\label{mono-hk}
h_{k}(n+1)\geq h_{k}(n).
\end{equation}
Furthermore, when $k\geq 2$ and $n\geq 0$,
%
\begin{equation}
\label{k^2-h-k}
k^{2}h_{k}(n)\leq n^{2}h_{k-1}(n).
\end{equation}
\end{lem}

\begin{proof} By definition, it is clear that for $k\geq 1$,
\begin{equation*}
1+\sum _{n=1}^{\infty }(g_{k}(n)-g_{k}(n-1))q^{n}=
\frac{1-q}{(q;q)_{k}^{2}}=\frac{1}{(q;q)_{k}(q^{2};q)_{k-1}},
\end{equation*}
which obviously has nonnegative power series coefficients. This yields  {\eqref{mono-gk}}.

Similarly, by  {\eqref{def-hmn}}, we see that for $k\geq 1$,
\begin{equation*}
\sum _{n=1}^{\infty }(h_{k}(n)-h_{k}(n-1))q^{n}=
\frac{(1-q)q^{2k}}{(q;q)_{k} (q^{2};q)_{k-1}}=
\frac{q^{2k}}{(q^{2};q)^{2}_{k-1}},
\end{equation*}
which also has nonnegative power series coefficients. Hence  {\eqref{mono-hk}} is valid.

We next prove  {\eqref{k^2-h-k}}. By  {\eqref{rec-gmn}}, we see that when
$k\geq 2$,
\begin{equation*}
k^{2} h_{k}(n)=k^{2}\sum _{i=0}^{\left \lfloor \frac{n}{k}\right
\rfloor -2}(i+1)h_{k-1}(n-ki-2).
\end{equation*}
From  {\eqref{mono-hk}}, we find that when $k\geq 2$ and $n\geq 0$,
\begin{align*}
k^{2} h_{k}(n)&\leq k^{2}\sum _{i=0}^{\left \lfloor \frac{n}{k}
\right \rfloor -2}(i+1)h_{k-1}(n)
\\[3pt]
&\leq k^{2}
\left (\left \lfloor \frac{n}{k}\right \rfloor -1\right )^2 h_{k-1}(n)
\\[3pt]
&\leq n^{2}h_{k-1}(n),
\end{align*}
as desired. This completes the proof.\end{proof}

We are now in a position to prove  {Theorem~\ref{lem-indu-im}}.

\begin{proof}[Proof of  {Theorem~\ref{lem-indu-im}}] (1) When $k=1$, it follows immediately
from  {\eqref{rec-fmn}} that $g_{1}(n)=n+1$ by noting that $g_{0}(0)=1$ and
$g_{0}(n)=0$ for $n\geq 1$. On the other hand, by the definition of
$h_{k}(n)$, we see that $h_{1}(0)=h_{1}(1)=0$, and $h_{1}(n)=1$ for
$n\geq 2$. Hence $g_{1}(n)\geq 21 h_{1}(n)$ when $n\geq 20$.

(2) When $k=2$, we first claim that when $n\geq 0$,
%
\begin{equation}
\label{ine-es-g2}
g_{2}(n)\geq \frac{n^{3}}{24}.
\end{equation}
Set $n=2t+j$, where $t\geq 0$ and $j=0$ or $1$. Notice that
$g_{1}(n)=n+1$, by  {\eqref{rec-fmn}}, we see that for $n\geq 0$,
\begin{eqnarray*}
g_{2}(n)&=&\sum _{i=0}^{\left \lfloor n/2\right \rfloor } (i+1)g_{1}(n-2i)
\nonumber
\\[3pt]
&=&\sum _{i=0}^{\left \lfloor n/2\right \rfloor } (i+1)(n-2i+1)
\nonumber
\\[3pt]
&\geq &\sum _{i=0}^{t} (i+1)(2t-2i+1)
\nonumber
\\[3pt]
&=&\frac{t^{3}}{3}+\frac{3t^{2}}{2}+\frac{13t}{6}+1.
\end{eqnarray*}
Hence, we derive that
%
\begin{equation}
\label{ine-est-f2n}
g_{2}(n)\geq \frac{t^{3}}{3}+\frac{3t^{2}}{2}+\frac{13t}{6}+1\geq
\frac{(t+1)^{3}}{3}.
\end{equation}
Since $n\leq 2t+2$, we deduce from  {\eqref{ine-est-f2n}} that for
$n\geq 0$,
\begin{equation*}
g_{2}(n)\geq \frac{(t+1)^{3}}{3}\geq \frac{n^{3}}{24}.
\end{equation*}
This yields  {\eqref{ine-es-g2}}.

On the other hand, since $h_{1}(n)=0$ or $1$ for $n\geq 0$, and by  {\eqref{k^2-h-k}}, we find that for $n\geq 0$,
%
\begin{equation}
\label{ine-h2-n24}
h_{2}(n) \leq \frac{n^{2}}{4}.
\end{equation}
Together with  {\eqref{ine-es-g2}}, we derive that for $n\geq 126$,
\begin{equation*}
g_{2}(n)\geq \frac{n^{3}}{24}\geq \frac{21n^{2}}{4}\geq 21 h_{2}(n).
\end{equation*}
Moreover, it can be checked that $g_{2}(n)\geq 21 h_{2}(n)$ when
$51\leq n\leq 125$. Hence we conclude that $g_{2}(n)\geq 21 h_{2}(n)$ for
$n\geq 51$.

(3) By suitable modification to the proof of  {\eqref{ine-es-g2}}, we can
show that
%
\begin{equation}
\label{ine-g3-43}
g_{3}(n)\ge \frac{n^{5}}{4320}.
\end{equation}

On the other hand, combining  {\eqref{k^2-h-k}} and  {\eqref{ine-h2-n24}}, we
find that when $n\geq 0$,
%
\begin{equation}
\label{ine-h3-n24}
h_{3}(n) \leq \frac{n^{4}}{36}.
\end{equation}
Hence by  {\eqref{ine-g3-43}} and  {\eqref{ine-h3-n24}}, we derive that when
$n\geq 2520$,
\begin{equation*}
g_{3}(n)\geq \frac{n^{5}}{4320}\geq \frac{21n^{4}}{36}\geq 21 h_{3}(n).
\end{equation*}
Furthermore, it is easy to check that $g_{3}(n)\geq 21h_{3}(n)$ when
$67\leq n\leq 2519$. Hence we conclude that
$g_{3}(n)\geq 21h_{3}(n)$ when $n\geq 67$.

(4) For $k\geq 4$, we will prove that $g_{k}(n)\geq 21h_{k}(n)$ when
$n\geq 0$ by induction on $k$.

When $k=4$, using the same method as above and after some tedious but straightforward
calculation, we deduce that when $n\geq 8$,
%
\begin{equation}
g_{4}(n)\geq \frac{1}{2903040}n^{7}.
\nonumber
\end{equation}
Here we omit the detail.

On the other hand, from  {\eqref{k^2-h-k}} and  {\eqref{ine-h3-n24}}, we deduce
that when $n\geq 0$,
%
\begin{equation}
h_{4}(n) \leq \frac{n^{6}}{576}.
\nonumber
\end{equation}
Hence, when $n\geq 105840$,
\begin{equation}
g_{4}(n)\geq \frac{n^{7}}{2903040}\geq \frac{21n^{6}}{576}\geq 21h_{4}(n).
\nonumber
\end{equation}
Furthermore, it can be checked that $g_{4}(n)\geq 21 h_{4}(n)$ when
$0\leq n\leq 105839$, so $g_{4}(n)\geq 21 h_{4}(n)$ for $n\geq 0$.

We now assume that there exists $k\geq 5$ such that
$g_{k-1}(n)\geq 21h_{k-1}(n)$ for $n\geq 0$. We aim to show that for
$n\geq 0$,
\begin{equation*}
g_{k}(n)\geq 21 h_{k}(n).
\end{equation*}
From  {\eqref{rec-fmn}} and  {\eqref{mono-gk}}, we derive that
\begin{eqnarray}
g_{k}(n)&=&\sum _{i=0}^{\left \lfloor \frac{n}{k}\right \rfloor }(i+1)g_{k-1}(n-ki)
\nonumber
\\[2pt]
&\geq & \sum _{i=0}^{\left \lfloor \frac{n}{k}\right \rfloor }(i+1)g_{k-1}(n-ki-2).
\nonumber
\end{eqnarray}
By the induction hypothesis, we have
\begin{equation*}
g_{k-1}(n-ki-2)\geq 21 h_{k-1}(n-ki-2).
\end{equation*}
Hence
%
\begin{equation}
g_{k}(n)\geq 21\sum _{i=0}^{\left \lfloor \frac{n}{k}\right \rfloor }(i+1)h_{k-1}(n-ki-2).
\nonumber
\end{equation}
From  {\eqref{rec-gmn}}, we have
%
\begin{equation}
h_{k}(n)=\sum _{i=0}^{\left \lfloor \frac{n}{k}\right \rfloor -2}(i+1)h_{k-1}(n-ki-2).
\nonumber
\end{equation}
It follows that for $n\geq 0$,
\begin{equation*}
g_{k}(n)\geq 21h_{k}(n).
\end{equation*}
This completes the proof.\end{proof}

We conclude this section with a proof of  {Conjecture~\ref{conj-opst3}}. Let
$N(\leq m,n)$ denote the number of partitions of $n$ with rank less than
or equal to $m$, and let $M(\leq m,n)$ denote the number of partitions
of $n$ with crank less than or equal to $m$. Bringmann and Mahlburg
\cite{{Bringmann-Mahlburg-2009}} conjectured that for $n\geq 1$ and
$m\leq 0$,
%
\begin{equation}
\label{inebre-1}
M(\leq m,n)\leq N(\leq m+1,n),
\end{equation}
which has been proved by Chen, Ji and Zang in
\cite{Chen-Ji-Zang-2017}.

By using the following two symmetries of ranks and cranks (see
\cite{Dyson-1969,Dyson-1989}):
%
\begin{equation}
\label{sym-cran}
M(m,n)=M(-m,n),
\end{equation}
and
%
\begin{equation}
\label{sym-ran}
N(m,n)=N(-m,n),
\end{equation}
it is not difficult to derive from  {\eqref{inebre-1}} that for
$n\geq 1$ and $m\geq 0$,
%
\begin{equation}
\label{inebre}
N(\leq m-1,n)\leq M(\leq m,n).
\end{equation}

We are now in a position to give a proof of  {Conjecture~\ref{conj-opst3}} by means of {Theorem~\ref{main-thm-n}} and  {Theorem~\ref{lem-indu}} as well as  {\eqref{chan-mao-upbound}}, {\eqref{inebre-1}},  {\eqref{sym-cran}} and  {\eqref{inebre}}.

\begin{proof}[Proof of  {Conjecture~\ref{conj-opst3}}] Setting $m=-2$ in  {\eqref{inebre-1}}, we see that for $n\geq 1$,
%
\begin{equation}
\label{ine--1}
M(\leq -2,n)\leq N(\leq -1,n).
\end{equation}
Setting $m=2$ in  {\eqref{inebre}}, we see that for $n\geq 1$,
%
\begin{equation}
\label{ine-m2}
M(\leq 2,n)\geq N(\leq 1,n).
\end{equation}
Subtracting  {\eqref{ine--1}} from  {\eqref{ine-m2}} leads to
%
\begin{equation}
\label{ineq-rancran}
N(0,n)+N(1,n) \leq M(-1,n)+M(0,n)+M(1,n)+M(2,n),
\end{equation}
for $n\geq 1$. By  {\eqref{sym-cran}}, we see that for $n\geq 1$,
\begin{equation*}
M(-1,n)=M(1,n),
\end{equation*}
and so  {\eqref{ineq-rancran}} becomes
%
\begin{equation}
\label{ineq-rancran-temp}
N(0,n)+N(1,n) \leq M(0,n)+2M(1,n)+M(2,n).
\end{equation}
By  {Theorem~\ref{main-thm-n}}, we see that for $n\geq 44$,
\begin{equation*}
M(0,n)\geq M(1,n)\geq M(2,n),
\end{equation*}
and so by  {\eqref{ineq-rancran-temp}}, we derive that for $n\geq 44$,
%
\begin{equation}
\label{ineq-rancran-r}
N(0,n)+N(1,n)\leq 4M(0,n).
\end{equation}
From  {\eqref{chan-mao-upbound}}, we see that for $n\geq 7$,
%
\begin{equation}
\label{chan-mao-u}
\mathop{\mathrm{ospt}}\nolimits (n)<\frac{p(n)}{4}+\frac{N(0,n)}{2}-
\frac{M(0,n)}{4}+\frac{N(1,n)}{2}.
\end{equation}
Applying  {\eqref{ineq-rancran-r}} in  {\eqref{chan-mao-u}}, we are led to
%
\begin{equation}
\label{chan-mao-u2}
\mathop{\mathrm{ospt}}\nolimits (n)<\frac{p(n)}{4}+
\frac{7M(0,n)}{4}
\end{equation}
for $n\geq 44$. Appealing to  {Theorem~\ref{lem-indu}}, we see that for
$n\geq 39$,
%
\begin{equation}
\label{ine-main-3}
21M(0,n)\leq p(n).
\end{equation}
Hence we arrive at
%
\begin{equation}
\label{chan-mao-u-f}
\mathop{\mathrm{ospt}}\nolimits (n)<\frac{p(n)}{3}
\end{equation}
for $n\geq 44$. Furthermore, it is easy to check that
$\mathop{\mathrm{ospt}}\nolimits (n)<\frac{p(n)}{3}$ when
$10\leq n\leq 43$. Thus, we complete the proof of  {Conjecture~\ref{conj-opst3}}.\end{proof}

\section{Conjectures}
\label{sec11}

Recall that a sequence $\{a_{i}\}_{1\le i\le n}$ is called log-concave
if for $2\le i\le n-1$, $a_{i}$ satisfies the following inequality:
\begin{equation*}
a_{i}^{2}\ge a_{i-1}a_{i+1}.
\end{equation*}
It is well known that if a sequence $\{a_{i}\}$ of positive integers is
log-concave, then $\{a_{i}\}$ is unimodal, see
\cite[P.124, Ex.50]{Stanley-1997}.

An interesting phenomenon occurs when we consider the log-concavity of
$M(m,n)$. In particular, for $72\le n\le 10000$ and
$72-n\le m\le n-72$ (tested with Mathematica), the following inequality
holds,
\begin{equation*}
M(m,n)^{2}\ge M(m-1,n)M(m+1,n).
\end{equation*}

In this case, we would like to make the following conjecture.

\begin{conj}%
\label{conj-m-log-concave}
For $n\geq 72$ and $72-n\le m\le n-72$,
%
\begin{equation}
\label{ine-log-mmn}
M(m,n)^{2}\ge M(m-1,n)M(m+1,n).
\end{equation}
In other words, for $n\geq 72$, the sequence
$\{M(m,n)\}_{|m|\le n-71}$ is log-concave.
\end{conj}
Obviously, this conjecture implies the sequence
$\{M(m,n)\}_{|m|\leq n-71}$ is unimodal when $n\geq 72$. More precisely,
when $n\ge 72$ and $1\le m\le n-71$,
%
\begin{equation}
\label{ineq-crank}
M(m-1,n)\ge M(m,n).
\end{equation}
It should be noted that the inequality  {\eqref{ineq-crank}} is also valid
when $n\ge 72$ and $n-70\le m\le n-1$, which follows immediately from the
following lemma. Hence if  {Conjecture~\ref{conj-m-log-concave}} is proved
to be true, then we could derive that the sequence
$\{M(m,n)\}_{|m|\leq n-1}$ is unimodal when $n\geq 44$ which is  {Corollary~\ref{unimod-cor}}.

\begin{lem}%
\label{log-uni}
For  $n\ge 2i\ge 4$,
%
\begin{equation}
\label{log-uni-equa}
M(n-i,n)=p_{i}(i),
\end{equation}
where $p_{r}(n)$ counts the number of partitions of $n$ with parts taken
from $\{2,3,\ldots ,r\}$, as defined in  {\eqref{gf-pkn}}.
\end{lem}
 \begin{proof} Let
$\lambda =(\lambda _{1},\lambda _{2},\ldots ,\lambda _{\ell })$ be a partition
of $n$ counted by $M(n-i,n)$. Let $n_{1}(\lambda )$ denote the number of
$1$'s in $\lambda $. We claim that $n_{1}(\lambda )=0$ when
$n\geq 2i$. Otherwise, by the definition of crank, we see that there are
exactly $n-i+n_{1}(\lambda )$ parts in $\lambda $ strictly larger than
$n_{1}(\lambda )$. Note that $n_{1}(\lambda )\geq 1$, so there are at least
$n-i+n_{1}(\lambda )$ parts in $\lambda $ not less than $2$. This leads
to
\begin{equation*}
n=\sum _{i=1}^{\ell }\lambda _{i} \ge 2(n-i+n_{1}(\lambda ))\ge 2n-2i+2
\ge n+2,
\end{equation*}
a contradiction. Hence $n_{1}(\lambda )=0$ and by the definition of crank,
we see that $\lambda _{1}=n-i$, which implies
\begin{equation*}
\sum _{i=2}^{\ell }\lambda _{i}= n-\lambda _{1}=i.
\end{equation*}
Moreover, $n_{1}(\lambda )=0$ implies that $\lambda _{\ell }\ge 2$. Hence
$(\lambda _{2},\ldots ,\lambda _{\ell })$ is counted by $p_{i}(i)$. Conversely,
let $\mu =(\mu _{1},\ldots ,\mu _{j})$ be a partition counted by
$p_{i}(i)$, notice that $n-i\ge i\ge \mu _{1}$, so the partition
$(n-i,\mu _{1},\ldots ,\mu _{j})$ is counted by $M(n-i,n)$. Hence the equality {\eqref{log-uni-equa}} holds. \end{proof}

We proceed to illustrate the inequality  {\eqref{ineq-crank}} is true when
$n\ge 72$ and $n-70\le m\le n-1$ with the aid of  {Lemma~\ref{log-uni}}. By
a direct calculation, it is easy to check that when $2\le i\le 70$,
\begin{equation*}
p_{i+1}(i+1)\ge p_{i}(i).
\end{equation*}
Thus from  {Lemma~\ref{log-uni}}, we see that when $n\ge 142$ and
$2\le i\le 70$,
\begin{equation*}
M(n-i-1,n)\geq M(n-i,n).
\end{equation*}
This leads to $M(m-1,n)\ge M(m,n)$ for $n\ge 142$ and
$n-70\le m\le n-2$. Note that $M(n-2,n)=1>0=M(n-1,n)$. Hence
$M(m-1,n)\ge M(m,n)$ is valid for $n\ge 142$ and $n-70\le m\le n-1$. After
checking the small cases for $72\le n\le 141$, we have verified
$M(m-1,n)\ge M(m,n)$ for $n\ge 72$ and $n-70\le m\le n-1$.

The similar phenomenon also occurs for $N(m,n)$. We have the following
conjecture.

\begin{conj}%
\label{ine-log-rank}
For $n\geq 73$ and $73-n\le m\le n-73$,
%
\begin{equation}
\label{ine-log-mmnN}
N(m,n)^{2}\ge N(m-1,n)N(m+1,n).
\end{equation}
In other words, for $n\geq 73$, the sequence
$\{N(m,n)\}_{|m|\le n-72}$ is log-concave.
\end{conj}

In \cite{Chan-Mao-2014}, Chan and Mao raised a problem of finding the condition
such that the inequality  {\eqref{ine-nm01-nmn}} holds. By the calculation
with Mathematica, we have the following conjecture.
%
\begin{conj}%
\label{conj-uni-rank}
For $n\ge 39$ and $1\le m\le n-2$,
%
\begin{equation}
\label{ine-nm01-nmn}
N(m-1,n)\ge N(m,n).
\end{equation}
\end{conj}
By the symmetry $N(m,n)=N(-m,n)$, we see that  {Conjecture~\ref{conj-uni-rank}} implies the sequence $\{N(m,n)\}_{|m|\le n-2}$ is unimodal
for $n\ge 39$.

It is clear that  {Conjecture~\ref{ine-log-rank}} implies the inequality  {\eqref{ine-nm01-nmn}} holds when $n\ge 73$ and $1\le m\le n-72$. We will
prove that the inequality  {\eqref{ine-nm01-nmn}} also holds when
$n\ge 73$ and $n-71\le m\le n-2$ by using the following lemma. Therefore,
the inequality  {\eqref{ine-nm01-nmn}} holds when $n\ge 73$ and
$1\le m\le n-2$. Furthermore, it is easy to check that  {\eqref{ine-nm01-nmn}} holds when $39 \leq n \leq 72$ and
$1\le m\le n-2$. Hence we could say that  {Conjecture~\ref{ine-log-rank}} implies {Conjecture~\ref{conj-uni-rank}}.

\begin{lem}%
\label{lem-exp-nn-i}
For  $n\ge 2i\ge 4$,
%
\begin{equation}
\label{lem-exp-nn-i-equal}
N(n-i,n)=\sum _{\ell =2}^{\lfloor (i+1)/2\rfloor }p(i-2\ell +1,\ell -1){.}
\end{equation}
where $p(n,r)$ denotes the number of partitions of $n$ with at most
$r$ parts, as defined in  {\eqref{equ-gen-pnkqn}}.
\end{lem}

\begin{proof} Let $\lambda =(\lambda _{1},\ldots ,\lambda _{\ell })$ be a
partition counted by $N(n-i,n)$. We claim that $\ell \geq 2$. Otherwise
if $\ell =1$, then $\lambda =(n)$ which is counted by $N(n-1,n)$, which
contradicts to the fact that $\lambda $ is counted by $N(n-i,n)$ where
$i\geq 2$. By the definition of rank, we see that
$\lambda _{1}=n-i+\ell $, so
\begin{equation*}
\sum _{j=2}^{\ell } (\lambda _{j}-1)=n-\lambda _{1}-\ell +1=i-2\ell +1
\geq 0.
\end{equation*}
Hence
$\mu =(\lambda _{2}-1,\lambda _{3}-1, \ldots ,\lambda _{\ell }-1)$ is a
partition counted by $p(i-2\ell +1,\ell -1)$. Furthermore,
$\ell \le (i+1)/2$.

Conversely, for $i\geq 2$ and
$2\leq \ell \leq \lfloor (i+1)/2\rfloor $, and let
$\mu =(\mu _{1},\mu _{2},\ldots ,\mu _{\ell -1})$ be a partition counted
by $p(i-2\ell +1,\ell -1)$, where $\mu _{\ell -1}\ge 0$. Note that
$n\geq 2i$, so $\mu _{1}\le i-2\ell +1\le n-i+\ell -1$. Hence
$\lambda =(n-i+\ell ,\mu _{1}+1,\mu _{2}+1,\ldots ,\mu _{\ell -1}+1)$ is
a partition counted by $N(n-i,n)$. Thus we arrive at  {\eqref{lem-exp-nn-i-equal}}. \end{proof}

We proceed to show that the inequality  {\eqref{ine-nm01-nmn}} holds when
$n\ge 73$ and $n-71\le m\le n-2$ with the aid of  {Lemma~\ref{lem-exp-nn-i}}. By a straightforward calculation, it is easy to check
that when $2\le i\le 71$,
\begin{equation*}
\sum _{\ell =2}^{\lfloor (i+2)/2\rfloor }p(i-2\ell +2,\ell -1)\ge
\sum _{\ell =2}^{\lfloor (i+1)/2\rfloor }p(i-2\ell +1,\ell -1).
\end{equation*}
Thus from  {Lemma~\ref{lem-exp-nn-i}}, we see that when $n\ge 144$ and
$2\le i\le 71$,
%
\begin{equation}
\label{ine-mn-i-n-le}
N(n-i-1,n)\geq N(n-i,n).
\end{equation}
This leads to $N(m-1,n)\ge N(m,n)$ for $n\ge 144$ and
$n-71\le m\le n-2$. It is easy to check that $N(m-1,n)\ge N(m,n)$ for
$73\le n\le 143$ and $n-71\le m\le n-2$. Hence $N(m-1,n)\ge N(m,n)$ when
$n\ge 73$ and $n-71\le m\le n-2$.

It should be noted that  {Conjecture~\ref{conj-m-log-concave}} and  {Conjecture~\ref{ine-log-rank}} have also been raised by Bringmann, Jennings-Shaffer
and Mahlburg \cite[Conjecture 4.3]{Bringmann-Jennings-Mahlburg-2021}.

 \vskip 0.2cm
\noindent{\bf Acknowledgments.} This work
was supported by   the National Science Foundation of China. We are greatly indebted to referees for their helpful suggestions that improved the presentation of this paper.

\end{document}